\documentclass[a4paper]{article}
\usepackage{amssymb}
\usepackage{amsthm}
\usepackage{amsmath}

\usepackage{mathrsfs}
\usepackage{dsfont}
\usepackage{xcolor}
\usepackage[T1]{fontenc}
\usepackage[utf8]{inputenc}
\usepackage{hyperref}
\newcommand\R{\mathbb{R}}

\newcommand\N{\mathbb{N}}

\newtheorem{thm}{Theorem}[section]
\newtheorem{cor}[thm]{Corollary}
\newtheorem{defi}[thm]{Definition}

\newtheorem{examp}[thm]{Example}
\newtheorem{lem}[thm]{Lemma}
\newtheorem{com}[thm]{Comment}
\newtheorem{prop}[thm]{Proposition}
\newtheorem{rk}[thm]{Remark}

\title{SET-VALUED STOCHASTIC INTEGRALS IN UMD SPACES AND APPLICATIONS}

\author{E. H.  Essaky \and M. Hassani \and C. E. Rhazlane \and \\
	Cadi Ayyad University\\ Poly-disciplinary Faculty of Safi\\
	Modeling and Combinatorics Laboratory\\
	Department of Mathematics and Computer Science\\ P.O. Box 4162, Safi, Morocco.\\ \\\hspace*{-0.7cm} e-mails : essaky@uca.ac.ma\qquad m.hassani@uca.ma\qquad charafeddine.rhazlane@gmail.com}
\date{}

\begin{document}

\maketitle

\begin{abstract}
The purpose of this paper is to study certain set-valued integrals in UMD Banach spaces and provide a compatible form of the martingale representation theorem for set-valued martingales. Under specific conditions, these martingales can be expressed using revised set-valued stochastic integrals with respect to a real standard Brownian motion $W = (W_t)_{t\in[0,T ]}$. 
Moreover, we prove the existence of solutions to the following set-valued backward stochastic differential equation of the form
\begin{eqnarray*}
	Y_t=\left(\xi+\int_t^TH_u du+\int_{[0,t]}^{\mathscr{R}}{Z}\,  dW_u\right) \circleddash \int_{[0,T]}^{\mathscr{R}}{Z}\, dW_u\quad a.s.,\quad t\in [0,T],
\end{eqnarray*}
where the right-hand side, of this equation, represents the Hukuhara difference of two quantities containing revised set-valued stochastic integrals, $\xi$ is a terminal set-valued function condition and $H$ is a set-valued function satisfying some suitable conditions.
\end{abstract}

\section{Introduction}

\indent Due to their connections with statistical and economic problems, the integration of set-valued functions has attracted increasing interest over time. Several approaches to defining the integration of set-valued functions can be distinguished. Firstly, Aumann \cite{Aumann} defined a set-valued integral using the set of all integrable selections. Hukuhara \cite{hukuhara} studied the integration of measurable functions whose values are compact convex sets, while Debreu \cite{debreu} employed an embedding method to consider the Bochner integral in an embedded Banach space. \\Let us mention here that other properties of set-valued integrals have been studied, notably in the work of Hiai and Umegaki \cite{Hiai}, where new concepts were introduced and played a key role in constructing the multivalued conditional expectation with respect to a $\sigma$-finite measure. Later, the notion of set-valued functions was extended to the stochastic case when Kisielewicz, in \cite{kis1997}, presented concepts of set-valued functional stochastic integrals with a selection property that is considered fundamental for stochastic inclusions.

This class of stochastic integrals does not yield set-valued stochastic processes nor exhibit properties analogous to those of point-valued stochastic integrals. Subsequently, Jung and Kim, in \cite{Jung}, proposed a definition of a set-valued stochastic integral with respect to a standard Brownian motion, building on Theorem $3.8$ of [\cite{papa}, Chapter 2]. Nevertheless, these set-valued stochastic integrals are not integrably bounded in the general case (see \cite{Michta}). To overcome this limitation, Kisielewicz introduced a new concept called the generalized set-valued stochastic integral in \cite{kis2014}.

Let $(\Omega,\mathscr{F},\mathbb{F},\mathbb{P})$ be a filtered probability space. Before defining our set-valued integrals in UMD Banach spaces, let us first clarify the definitions of the following set-valued integrals:
\begin{itemize}
\item For a given measurable, adapted, and $2$-integrably bounded set-valued function $H:[0,T]\times \Omega \rightarrow\mathscr{K}(\R^m)$, the Lebesgue set-valued stochastic integral of $H$ over $[s, t]$, denoted as $(\ell_{\mathbb{F}})\int_s^t H_udu$, is defined in [\cite{set-stochastic-ch}, Chapter $4$, Section $3$] as a set-valued random variable satisfying
\begin{eqnarray}\label{lebint}
S_{\mathscr{F}_t}^2\left((\ell_{\mathbb{F}})\int_s^t H_udu\right)=\overline{dec}_{\mathscr{F}_t}\left(\left\lbrace \int_s^tg(u,.)du:g\in S_{\Sigma_{\mathbb{F}}}^2(H)\right\rbrace\right),
\end{eqnarray}
where $\Sigma_{\mathbb{F}}$ is the $\sigma$-algebra of all measurable and adapted subsets of\; $[0,T]\times \Omega$, $S_{\mathscr{F}_t}^2$ refers to the subtrajectory integrals described in Definition \ref{def},  $\overline{dec}_{\mathscr{F}_t}$ denotes the closed decomposable hull as defined in Definition \ref{def2}, and
\begin{eqnarray*}\label{sigma}
S_{\Sigma_{\mathbb{F}}}^2(H)= \Big\{  g\in L^2([0,T]\times \Omega,\Sigma_{\mathbb{F}};\R^m):g(v,\omega)\in H(v,\omega)\;a.e.(v,\omega)\Big\}.
\end{eqnarray*}
For further reading on the class of set-valued stochastic integrals with respect to the Lebesgue measure, see [\cite{kis2012}, Section 3].
\item For a nonempty subset $\tilde{\mathscr{G}}$ of the space $L^2([0,T]\times \Omega,\Sigma_{\mathbb{F}};\R^m)$, the generalized set-valued stochastic integral of $\tilde{\mathscr{G}}$ over $[0,t]$ is understood by \cite{kis2014} as an $\mathscr{F}_t$-measurable set-valued random variable satisfying
\begin{eqnarray}\label{geneint}
S_{\mathscr{F}_t}^2\left(\int_0^t \tilde{\mathscr{G}} dW_u\right)=\overline{dec}_{\mathscr{F}_t}\left(\left\lbrace \int_0^tg_udW_u:g\in \tilde{\mathscr{G}}\right\rbrace\right),
\end{eqnarray}
where $W=(W_t)_{t\in [0,T]}$ is a standard Brownian motion defined on $(\Omega,\mathscr{F},\mathbb{F},\mathbb{P})$.
\end{itemize}

In this work, our main purpose is to study certain set-valued integrals in UMD Banach spaces and provide a compatible form of the martingale representation theorem for set-valued martingales. To achieve this, we first need to choose an appropriate space that offers the right setting for establishing a robust theory of stochastic integration of adapted processes taking values in a UMD space  $E$. We have selected the Banach space $L^p(\Omega, \gamma(0, T; E))$ for this purpose. Additionally, we aim to address the challenges associated with this choice.

%for every $p\in (1,\infty)$, the Banach space $L^p(\Omega, \gamma(0, T; E))$ which provide a right setting to establish a rich theory of stochastic integration of adapted processes with values in a UMD space $E$ [see, for example, the works by \cite{neerven0, neerven3}].
In order to give a more comprehensive view, we consider the following backward stochastic evolution equation
\begin{eqnarray}\label{bdeq1}
	\begin{cases}
		dU_t= f(t)dt+V_tdW_t,\quad t\in [0,T],\\
		U_T=\mathfrak{T},
	\end{cases}
\end{eqnarray}
where the existence of solutions has been proved, in Lü and Neerven \cite{BSEE.B}, by virtue of the martingale representation theorem. It should be noted that the integral with respect to time $t$ in (\ref{bdeq1}) is generally done by abuse of notation, and for the sake of readability it has been used. Starting from this point-valued case,  natural questions then arise:
% We then found ourselves faced with the resolution of situations related to the framework chosen for our study:
\begin{itemize}
\item   How to suitably define the integrals of set-valued stochastic processes, and what properties do they have?
\item How to define set-valued stochastic integrals with respect to the Brownian motion $W$ in the framework of UMD spaces? What properties do they have? Can they represent the set-valued martingale?
\end{itemize}
%As far as the theory of stochastic in UMD spaces allows us, our project is interested in providing and studying other set-valued integrals. 
To address the above questions, we consider the space $\mathcal{L}_{\mathbb{F}}^p([0,T]\times \Omega;E)$ of all $E$-valued measurable and adapted mappings defining elements of $L^p(\Omega,\gamma(0,T;E))$, where $\gamma(0,T;E)$ is the space of $\gamma$-radonifying operators form $L^2([0,T])$ into $E$. We then establish a concept concerning the $\gamma$-\;set-valued stochastic integral of a set-valued function $H:[0,T]\times \Omega\longrightarrow \mathscr{K}(E)$ by introducing the following set
\begin{eqnarray*}
\Lambda^{0,T}_{p,H}= \Big\{ f\in\mathcal{L}_{\mathbb{F}}^p([0,T]\times \Omega;E):f(v,\omega)\in H(v,\omega)\;a.e.(v,\omega)\in [0,T]\times \Omega\Big\}.
\end{eqnarray*}
This set is used to define the $\gamma$-\;set-valued stochastic integrals of $H$. Furthermore, by imposing additional conditions on the set $\Lambda^{0,T}_{p,H}$, we can derive further properties. In some cases, the subsets  $S_{\Sigma_{\mathbb{F}}}^2(H)$ described in Definition \ref{def} and $\Lambda^{0,T}_{p,H}$ will be comparable with respect to inclusion.

Concerning the integral with respect to a standard Brownian motion $W$, we can generalize the generalized set-valued stochastic integral to the framework of UMD spaces. We define the set-valued stochastic integral, $\int_0^t \mathscr{G}\, dW_u$, over $[0,t]$ of a nonempty subset $\mathscr{G}$ contained in $\mathcal{L}_{\mathbb{F}}^p([0,T]\times \Omega;E)$ by
\begin{eqnarray*}\label{g-s-v-int}
S_{\mathscr{F}_t}^p\left(\int_0^t \mathscr{G}\, dW_u\right)=\overline{dec}_{\mathscr{F}_t}\left(\left\lbrace \int_0^tg_udW_u:g\in \mathscr{G}\right\rbrace\right),
\end{eqnarray*}
where the closed decomposable hull is taken in the space $L^p(\Omega,\mathscr{F}_t;E)$. However, this set-valued stochastic integral was considered unsuitable when Zhang and Yano proved in \cite{J2022} that the martingale representation theorem for set-valued martingales \cite{kis2014} imposes such strong assumptions  that it rules out non-degenerate set-valued martingales, making it applicable only in the degenerate case.
Moreover, the revised representation theorem [\cite{J2022}, Theorem $2$] presented only equivalences highlighting the limitations of the generalized set-valued stochastic integrals and the conditions that a set-valued martingale must satisfy in order to have a representation analogous to the point-valued case. To resolve this situation, we extend these generalized set-valued integrals by considering nonempty subsets $\mathfrak{S}$ contained in the Cartesian product $E\times \mathcal{L}_{\mathbb{F}}^p([0,T]\times \Omega;E)$, and define the revised set-valued stochastic integrals over $[0,t]$ by
\begin{eqnarray*}
S_{\mathscr{F}_t}^p\left(\int_{[0,t]}^{\mathscr{R}} \mathfrak{S}\, dW_u\right)=\overline{dec}_{\mathscr{F}_t}\left(\left\lbrace x+\int_0^tg_udW_u:(x,g)\in \mathfrak{S}\right\rbrace\right).
\end{eqnarray*}
Under appropriate conditions, we opt to represent set-valued martingales using the revised set-valued integral, as this representation aligns with the set-valued context, adheres to the martingale representation theorem in classical stochastic analysis, and offers sufficient flexibility to extend (\ref{bdeq1}) and provide a solution to the following Set-Valued Backward Stochastic Differential Equation 
\begin{eqnarray}\label{sveq1}
Y_t=\left(\xi+\int_t^TH_u du+\int_{[0,t]}^{\mathscr{R}}{Z}\,  dW_u \right)\circleddash \int_{[0,T]}^{\mathscr{R}}{Z}\,  dW_u\quad a.s.\; \omega\in\Omega,\quad t\in [0,T],
\end{eqnarray}
where the right-hand side of (\ref{sveq1}) represents the Hukuhara difference of two quantities (see Definition \ref{defi201} below), $\xi$ is a set-valued terminal condition $\mathscr{F}_T$-measurable and $p$-integrably bounded, with $p$ is a fixed real number greater than $1$.
\\
\par Let us describe our plan. In section $2$, {we present most of the necessary concepts and terminology and prove some results that will be useful for the analysis required to demonstrate and illustrate some of the main outcomes. Section 3 explores the extension of certain properties from \cite{kis2013, kis2014} to infinite dimensions, using tools from vector integration and stochastic integration in UMD spaces.} Additionally, we derive representations related to the revised set-valued integrals, assuming the separability of the probability space.\\
In Section $4$, we introduce the $\gamma$-set-valued stochastic integrals for set-valued functions. We establish the Chasles relation for these integrals. We prove a representation using a sequence in the set $\Lambda_{p,H}^{0,T}$, and present additional results by equipping the space $E$ with geometric properties.
\\Section $5$ is devoted to the martingale representation theorem for given set-valued martingales in UMD spaces, along with the presentation of additional results.\\ 
In Section $6$, we establish the existence of a solution to the set-valued backward stochastic differential equation (\ref{sveq1}). We also examine the case where $E$ has martingale type $2$.

%%%%%%%%%%%%%%%%%%%%%%%%%%%%%%%%%%%%%Section2%%%%%%%%%%%%%%%%%%%%%%%%%%%%%%%%%%%%%%%%%%%%%%%
%we provide some  and ingredients that will be effective during the analysis in the sequel. Essentially,

\section{Preliminaries}\label{sec:2}
In this section, we will recall some notions related to Banach space analysis, introduce certain prerequisites associated with stochastic integration in UMD Banach spaces, and highlight some tools from multivalued analysis.

\subsection{Notations}
Throughout this paper, the normed vector spaces are assumed to be real, and we will always identify Hilbert spaces with their dual by means of the Riesz representation theorem.  \\\hspace*{0.2cm}

In the sequel to this paper, we will use the following notation:\\\hspace*{0.4cm}

$\bullet$ Let $(E, \left\Vert \cdot \right\Vert_E)$ be a separable Banach space. The symbol $0_E$ will denote the additive identity element of $E$. The Borel $\sigma$-field of $E$ equipped with the norm topology will be denoted by $\mathscr{B}(E)$. The dual of $E$ will be denoted by $E^*$, and the duality between $E$ and $E^*$ will be represented by $\langle x,x^*\rangle$. \\\hspace*{0.2cm}

% Let $(E,\left\Vert\,.\,\right\Vert_E)$ be a separable Banach space. $0_E$ will denote the additive identity element of $E$. $\mathscr{B}(E)$ will denote the Borel $\sigma$-field of $E$ equipped with the norm topology. The dual of $E$ will be denoted by $E^*$, and the duality between $E$ and $E^*$ will be represented by $\langle x,x^*\rangle$.
$\bullet$ Let $\mathscr{H}$ be a Hilbert space with inner product $\langle \cdot, \cdot \rangle_{\mathscr{H}}$. The set of all bounded linear operators from $\mathscr{H}$ into $E$ will be denoted by $\mathcal{L}(\mathscr{H}, E)$.  \\\hspace*{0.2cm}

%Let $\mathscr{H}$ be a Hilbert space with inner product $\langle .,.\rangle_{\mathscr{H}}$. 
%$\mathcal{L}(\mathscr{H},E)$ will denote the set of all bounded linear operators from $\mathscr{H}$ into $E$.
$\bullet$ The set of all nonempty closed (resp. nonempty closed convex) subsets of $E $ will be denoted by   $\mathscr{K}(E)$ and $\mathscr{K}_c(E)$, respectively,  while $\mathscr{K}_{cwcmpt}(E)$ will refer to the set of all nonempty convex weakly compact subsets of $E$.\\\hspace*{0.2cm}

%We shall denote by  $\mathscr{K}(E)\;\left(\textrm{resp.}\,\mathscr{K}_{c}\left( E\right)\right)$ the set of all nonempty closed ( resp. nonempty closed convex ) subsets of $E$. While we assign $\mathscr{K}_{cwcmpt}(E)$ to indicate the set of all nonempty convex weakly compact subsets of $E$.\\\hspace*{0.4cm}
$\bullet$ Let $A$ and $B$ two subsets of the space $\mathscr{K}\left( E\right)$. The Hausdorff distance $\varrho\left( A,\,B\right)$ is defined as follows:
\begin{eqnarray*}
	\varrho\left( A,\,B\right) :=\,\max\left(\tilde{\varrho}\left( A,\,B\right) ,\,\tilde{\varrho}\left( B,\,A\right)\right),
\end{eqnarray*}
where $$\tilde{\varrho}\left( A,\,B\right)=\displaystyle{\sup_{ x\in A}}\;\displaystyle{\inf_{ y\in B}}\left\Vert x-y\right\Vert_E,$$
and $\tilde{\varrho}\left( B,\,A\right)$ can be expressed by changing the roles of $A$ and $B$.\\
We will use the notation 
\begin{eqnarray}\label{n.hausdorff}
\left\vert A\right\vert_{\varrho}:=\varrho(A,\left\lbrace 0_E\right\rbrace)=\sup_{x\in A}\left\Vert x\right\Vert_E.	
\end{eqnarray}

$\bullet$ $\N$ will denote the set of all natural numbers,  $T$ and $p$ be two fixed real numbers such that $T>0$ and $p> 1$.\\\hspace*{0cm}

$\bullet$ For two positive real quantities $L_1$ and $L_2$, we shall write $L_1\lesssim_{p} L_2$ (resp. $L_1\lesssim_{p,E} L_2$) to express that there exists a positive constant $\kappa_p$ (resp. $\kappa_{p,E}$) depending only on $p$ (resp. on $p$ and $E$) such that 
$$L_1\leq \kappa_p L_2\,\, (\text{resp.}\,\, L_1\leq \kappa_{p,E} L_2).$$ The equivalence $L_1\eqsim_{p,E} L_2$ will mean that $L_1\lesssim_{p,E} L_2$ and $L_2\lesssim_{p,E} L_1$.\\\hspace*{0.2cm}

$\bullet$ $\left(\Omega,\mathscr{F},\mathbb{P},\mathbb{F}\right)$ will represent a filtered probability space where we set $\mathbb{F}=\Big(\mathscr{F}_t\big)_{T\geq t\geq 0},$ with $\mathscr{F}_T=\mathscr{F}$.\\\hspace*{0.2cm}

$\bullet$  $\pi_s,\,s\in [0,T],$ will denote the mapping that associates each $E$-valued process $(U_t)_{0\leq t\leq T}$ with $U_s$. \\\hspace*{0.2cm}

$\bullet$  $\Sigma_{\mathbb{F}}$ will denote the $\sigma$-algebra of all measurable and adapted subsets of\; $[0,T]\times \Omega$.
%%%%%%%%%%%%%%Section2-Subsection1%%%%%%%%%%%%%%%%%%%%%%%%%%
\subsection{Bochner and Pettis Integrability}
Let $(\mathcal{S},\mathscr{A},\mu)$ be a $\sigma-$finite measure space, and let $q$ be a real number such that $q\geq 1$.\\\hspace*{1cm}

\noindent We first note that there are various possible definitions of measurability for functions, as discussed in [\cite{neerven1}, Section 1.1]. In this paper, we consider functions that map to separable Banach spaces, and thus we focus on the classical notion of measurability. Specifically, a function $\mathfrak{h}: \mathcal{S} \rightarrow E$ is said to be measurable (or $\mathscr{A}$-measurable) if, for every $\mathcal{C} \in \mathscr{B}(E)$ the preimage $\mathfrak{h}^{-1}(\mathcal{C})$ is measurable.\\\hspace*{1cm}

We say that two functions are $\mu$-versions of each other if they agree $\mu$-almost everywhere.
%Let us call two functions which agree $\mu$-almost everywhere $\mu$-versions of
each other.
\begin{defi}
Let $f:\mathcal{S}\longrightarrow E$ be a $\mu$-version of a measurale function.
$f$ is called Bochner $q$-integrable if we have 
\begin{eqnarray*}
\left\Vert f\right\Vert^q_{L^q(\mathcal{S};E)}:=\int_{\mathcal{S}}\left\Vert f\right\Vert^q_Ed\mu<\infty,
\end{eqnarray*}
where the integral is in the Lebesgue sense.
\end{defi}
The linear space $L^q(\mathcal{S};E)$ consisting of (classes of a.e. equal) Bochner $q$-integrable functions $f:\mathcal{S}\longrightarrow E$, endowed with the norm $\left\Vert .\right\Vert_{L^q(\mathcal{S}; E)}$, is a Banach space.\\

For a given sub-$\sigma$- algebra $\tilde{\mathscr{A}}$ of $\mathscr{A}$, we will write $L^q(\mathcal{S},\tilde{\mathscr{A}};E)$ to express the closed linear subspace of $L^q(\mathcal{S};E)$ consisting of all equivalence classes of functions where each of them admits a representative that is $\mu_{\mid \tilde{\mathscr{A}}}$\,-version of a $\tilde{\mathscr{A}}$-measurable function.\\

We notice that the $\mu$-simple functions are dense in the Bochner space $L^q(\mathcal{S};E)$, and for any $\mu$-simple function $h=\displaystyle{\sum_{n=1}^N\mathds{1}_{A_n}x_n}$, we define the Bochner integral of $g$ as follows:
\begin{eqnarray}\label{formula3}
\int_{\mathcal{S}}^{Bochner}hd\mu:=\sum_{j=1}^{N}\mu(A_j)x_{j}.	
\end{eqnarray}
The expression in (\ref{formula3}) is independent of the representation of $h$.\\\hspace*{1cm}

\noindent The Bochner integral of a function $f \in L^1(\mathcal{S}; E)$ is defined as the limit of the Bochner integrals of $\mu$-simple functions $(f_n)_{n \geq 1}$, where $(f_n)_{n \geq 1}$ converges to $f$ in the norm of $L^1(\mathcal{S}; E)$. Additionally, for any $A \in \mathscr{A}$, the Bochner integral of $f$ over $A$ is simply the Bochner integral of $\mathds{1}_A f$.\\

Before introducing a more general integral, we first need to present the following notion.
\begin{defi}  A function $f:\mathcal{S}\rightarrow E$ is said to be weakly in $L^q(\mathcal{S})$ if the real function $\langle f,x^*\rangle$ belongs to the Lebesgue space $L^q(\mathcal{S})$ for every $x^*\in E^*$.
\end{defi}
Next, we will describe the notion of the Pettis integral.
%The following Definition-Theorem is due to Pettis \cite{pettis}.
\begin{defi}\label{def210} Let $f:\mathcal{S}\rightarrow E$ be a function weakly in $L^1(\mathcal{S})$. 
$f$ is called Pettis integrable if for every $A\in \mathscr{A}$, there exists $\nu_f(A)\in E$ such that for every $x^*\in E^*$
$$\langle \nu_f(A),x^*\rangle=\int_{A}\langle f,x^*\rangle d\mu.$$
In this case, the set function $\nu_f : \mathscr{A}\rightarrow E$ is called the Pettis integral of $f$ with respect to $\mu$, and
$\nu_f(A)$ is called the Pettis integral of $f$ over $A\in \mathscr{A}$ with respect to $\mu$, notation $\int_A^{Pettis}fd\mu:=\nu_f(A)$.
\end{defi}
It should be noted that the set of all Pettis integrable functions from $\mathcal{S}$ to 
$E$ forms a linear space. Moreover, the Pettis integral operator possesses the linearity property, and Bochner integrability implies Pettis integrability according to the identity [\cite{neerven1}, ($1.2$)].

\begin{rk}
Let $f: \mathcal{S} \rightarrow E$ be a $\mu$-version of a measurable function, assumed to be weakly in $L^2(\mathcal{S})$. According to Theorem $1.2.37$ in \cite{neerven1}, for any $g \in L^2(\mathcal{S})$, the function $s \mapsto g(s)f(s)$ is Pettis integrable. This allows us to define the linear operator
\begin{eqnarray}\label{map1}
	\begin{array}{ccccc}
		\mathscr{I}_f & : &  L^2(\mathcal{S})& \longrightarrow & E\\
		& & g & \longmapsto & \displaystyle	\mathscr{I}_f(g)= \int_{\mathcal{S}}^{Pettis}g(u)f(u)d\mu(u).
	\end{array}
\end{eqnarray}

\noindent Moreover, by the closed graph theorem, the operator $\mathscr{I}_f$ is bounded and $\mathscr{I}_f = \mathscr{I}_h$ for any function $h$ that is equal to $f$ almost everywhere.
\end{rk}
%We shall call $\mathscr{I}_f$ the Pettisintegral operator with kernel $f$.

%Next, let $f:[s,t]\longrightarrow E$ be a function strongly $\lambda$-measurable, weakly in $L^2([s,t])$, and let us consider the following useful mapping:
%\begin{eqnarray}\label{map1}
%\begin{array}{ccccc}
%\mathscr{I}_f & : &  L^2([s,t])& \longrightarrow & E\\
% & & g & \longmapsto & \displaystyle\int_{[s,t]}^{Pettis}g(u)f(u)d\lambda(u).
%\end{array}
%\end{eqnarray}
%Due to Theorem $1.2.37$ in \cite{neerven1}, $\mathscr{I}_f$ is well defined and linear. The closed graph theorem asserts that $\mathscr{I}_f$ is bounded. Moreover, $\mathscr{I}_f=\mathscr{I}_{h}$ for any function $h$ equal to $f$ almost everywhere.

\subsection{Functions Defining a $\gamma$-Radonifying Operator}

The $\gamma$-radonifying operators play a crucial role in developing the theory of Gaussian measures in Banach spaces as well as stochastic integrals on some Banach spaces (cf. \cite{kuo, dett, brze1, brze2, neerven0}), and they have been used to characterize some nice geometric structure, like type $2$ and cotype $2$ (see \cite{typeco}).\\

Let $(\gamma_n')_{n\geq 1}$ be a sequence of independent standard Gaussian random variables on a probability space $(\Omega',\mathscr{F}',\mathbb{P}')$. 
\begin{defi}
An operator in $\mathcal{L}(\mathscr{H}, E)$ is said to be of finite rank if it is a linear combination
of operators of the form $h\otimes e$, where $h\in \mathscr{H}$, $e\in E$, and $h\otimes e$ is defined by $$(h\otimes e)(h')=\langle h,h'\rangle_{\mathscr{H}}e\quad \textrm{ for every }h'\in \mathscr{H}.$$
The space of all finite rank operators contained in $\mathcal{L}(\mathscr{H}, E)$ is denoted by $\mathscr{H}\otimes E$.
\end{defi}
\vspace{0.4cm}
\noindent Due to the Gram-Schmidt orthogonalization argument, every finite rank operator can be represented in the form $\displaystyle{\sum_{j=1}^kh_j\otimes e_j}$, where $(h_j)_{1\leq j\leq k}$ is orthonormal in $\mathscr{H}$ and $(e_j)_{1\leq j\leq k}$ is in $E$.
\begin{defi}
The Banach space $\gamma(\mathscr{H}, E)$ of the $\gamma$-radonifying operators from $\mathscr{H}$ into $E$ is defined as the completion of $\mathscr{H}\otimes E$
with respect to the norm
\begin{eqnarray}\label{gamma1}
\left\Vert \sum_{j=1}^kh_j\otimes e_j\right\Vert_{\gamma(\mathscr{H},E)}:=\left\Vert \sum_{j=1}^k\gamma_j'e_j\right\Vert_{L^2(\Omega';E)},
\end{eqnarray}
where it is assumed that $h_1,..., h_k$ are orthonormal in $\mathscr{H}$.
\end{defi}

\noindent By virtue of Proposition $3.15$ in \cite{survey} we have for all $\mathbb{T}\in \gamma(\mathscr{H}, E)$
\begin{eqnarray}\label{gamma2}
\left\Vert \mathbb{T}\right\Vert_{\gamma(\mathscr{H},E)}:= \left(\sup \mathbb{E}'\left\Vert \sum_{j=1}^N \gamma'_j \mathbb{T} h_j\right\Vert_E^2\right)^{\frac{1}{2}},
\end{eqnarray}
where the supremum is taken over all finite orthonormal systems ${h_1,..., h_N}$ in $\mathscr{H}$.\\

\noindent We notice that the right-hand sides of (\ref{gamma1}) and (\ref{gamma2}) do not depend on $(\gamma_n')_{n\geq 1}$ and the probability space $(\Omega',\mathscr{F}',\mathbb{P}')$ according to a uniqueness property of Fourier
transforms [\cite{st.ev.eq}, Theorem $2.8$].\\ Moreover, the right-hand side of (\ref{gamma1}) is independent of the representation of the finite rank operator as long as we choose $(h_j)_{1\leq j\leq k}$ orthonormal in $\mathscr{H}$ (see [\cite{survey}, page $10$]).\\

Let us consider the measure space  $([s,t],\mathscr{B}([s,t]),\lambda)$ where $s,t$ are two real numbers such that $0\leq s<t\leq T$, and $\lambda$ is the Lebesgue measure on the $\sigma$-field $\mathscr{B}([s,t])$.\\

The spaces of $\gamma$-radonifying operators principally involved in our work are those for  $\mathscr{H}$ of the form $L^2([s,t])$, and we will use $\gamma(s,t;E)$ as a shorthand for $\gamma(L^2([s,t]),E)$.

%The spaces of $\gamma$-radonifying operators which principally involved in our work are for $\mathscr{H}$ of the form $L^2([s,t])$, and we will write $\gamma(s,t;E)$ instead of $\gamma(L^2([s,t]),E)$.

%%%%%%%%%%%%%%%%%%%%%%%%%%%%%%%%%%%%%Section2-Subsection3%%%%%%%%%%%%%%%%%%%%%%%%%%%%%%%%%%%

\begin{defi}\label{def241}
A function $f:[s,t]\longrightarrow E$ is said to define an element of $\gamma(s,t;E)$ if it satisfies the following conditions:
\begin{itemize}
\item[(i)]
$f$ is $\lambda$-version of a $\mathscr{B}([s,t])$-measurable function, and is weakly in $L^2([s,t])$;
\item[(ii)]
the Pettis integral operator $\mathscr{I}_f$ defined as in (\ref{map1}), is $\gamma$-radonifying.
\end{itemize}
In this case we write $f\in \gamma(s,t;E)$ and use the abbreviation $\left\Vert f\right\Vert_{\gamma(s,t;E)}:=\left\Vert \mathscr{I}_f\right\Vert_{\gamma(s,t;E)}$.
\end{defi}
\begin{rk}\label{rk241}
Let $f$ (resp. $h$) be a function belonging to $\gamma(s,t;E)$ (resp. $\gamma(0,T;E)$). Taking into account [\cite{neerven2}, Examples $9.1.11$-$9.1.12$], we have respectively 
$$
\left\Vert f\right\Vert_{\gamma(s,t;E)}=\left\Vert \mathds{1}_{[s,t]}f\right\Vert_{\gamma(0,T;E)} $$ 
and
$$
\left\Vert h_{\mid [s,t]}\right\Vert_{\gamma(s,t;E)}\leq\left\Vert h\right\Vert_{\gamma(0,T;E)},$$ 
where $h_{\mid [s,t]}$
is the restriction of $h$ on $[s,t]$.
\end{rk}

Next, we present the computation of the norm of a particular $\gamma$-radonifying operator.
\begin{examp}\label{examp1}
Let $h$ be a function belonging to $L^2([s,t])$ and $e$ be an element of the space $E$.
\begin{itemize}
\item Assume $\vert h(u)\vert>0$ for almost everywhere $u\in [s,t]$ and let's take $h_1=c_1h$ where $c_1:=1/\left\Vert h\right\Vert_{L^2([s,t])}$ is a normalising constant. We have
\begin{align*}
\left\Vert h\otimes e\right\Vert_{\gamma(s,t;E)}& =c_1^{-1}\left\Vert h_1\otimes e\right\Vert_{\gamma(s,t;E)}  \\ & =\left\Vert h\right\Vert_{L^2([s,t])}\left(\mathbb{E}'\left(\left\Vert \gamma'_1e\right\Vert_E^2\right)\right)^{\frac{1}{2}}=\left\Vert h\right\Vert_{L^2([s,t])}\left\Vert e\right\Vert_E.
\end{align*}
\item Assume that $h(u)=0$ for almost everywhere $u\in [s,t]$ and let us take $h_2$ of norm $1$ in $L^2([s,t])$. It follows that 
$$\left\Vert h\otimes e\right\Vert_{\gamma(s,t;E)}=\left\Vert h_2\otimes 0_E\right\Vert_{\gamma(s,t;E)}=\left(\mathbb{E}'\left(\left\Vert \gamma'_20_E\right\Vert_E^2\right)\right)^{\frac{1}{2}}=0.$$
\end{itemize}
Therefore, for every $(e,h)\in E\times L^2([s,t])$, we have $$\left\Vert h\otimes e\right\Vert_{\gamma(s,t;E)}=\left\Vert h\right\Vert_{L^2([s,t])}\left\Vert e\right\Vert_E.$$
\end{examp}

\par Now, let us consider the following mapping:
$$m_{s,t}:\phi\in L^2([0,T])\longmapsto \int_{[s,t]} \phi(u)du.$$
Obviously, every element $\tilde{\mathbb{T}}$ of $\gamma(\R,E)$ can be expressed as follows: $\tilde{\mathbb{T}}=1\otimes \tilde{\mathbb{T}}(1)$, where $1$ is the multiplicative identity element of $\R$. Thus, by virtue of (\ref{gamma2}) the linear operator $\mathbb{T}\in \gamma(\R,E)\longmapsto \mathbb{T}(1)$ establishes that $\gamma(\R,E)=E$ isometrically.
\\Next, according to Kalton-Weis extension theorem [\cite{neerven2}, Theorem $9.6.1$], the mapping $m_{s,t}\otimes I_E$, which associates any operator $\phi\otimes e$ with $(\int_{[s,t]}\phi(u)du) e$, has a unique extension to a bounded linear operator $\tilde{m}_{s,t}$ from $\gamma(0,T;E)$ to $E$, with the same norm $\sqrt{t-s}$ as $m_{s,t}$.

For the sake of readability, we will use the notation $\int _s^t f(u)du$,  with a slight abuse of notation, instead of the expression  $\tilde{m}_{s,t}(\mathscr{I}_f)$, for any function $f$ that defines an element of $\gamma(0,T;E)$.\\

In the following proposition, we list some useful properties satisfied by this integral notation. We only need to prove equality (\ref{ident212}). The linearity of the integral, expressed by (\ref{ident216}), and Chasles's relation (\ref{ident211}) for the integral are direct consequences of [\cite{BSEI}, Proposition 2.8].
\begin{prop}\label{prop245}
Let $\alpha$ be a real number, $s'$ be a real number with $s<s'<t$, and $f,g$ be two functions defining elements of $\gamma(0,T;E)$. The following properties hold:
\begin{eqnarray}
\int_s^t (f+\alpha g)(u)du &=& \int_s^t f(u)du+\alpha \int_s^t g(u)du,\label{ident216}\\
\int_s^t f(u)du &=& \int_s^{s'} f(u)du + \int_{s'}^t f(u)du,\label{ident211}\\
\int_s^t f(u)du &=& \int_0^{T} \mathds{1}_{[s,t]}f(u)du,\label{ident212}
\end{eqnarray}
where integrals of the functions are images of the associated Pettis integral operators by $\tilde{m}_{.,.}$.
\end{prop}

\begin{proof} 
Let $(\mathbb{T}_n)_{n\geq 1}$ be a sequence converging in $\gamma(0,T;E)$ to $\mathscr{I}_f$, such that for each $n\in \N$ $$\mathbb{T}_n:=\sum_{j=1}^{k_n}h_j^{n}\otimes e_j^n,$$ where $k_n\in \N,\,h_j^{n}\in L^2([0,T]),$ and $e_j^n\in E$. Next, we consider for each $n\in\N$ $$\tilde{\mathbb{T}}_n:=\displaystyle\sum_{j=1}^{k_n}\left(\mathds{1}_{[s,t]}h_j^{n}\right)\otimes e_j^n.$$
By virtue of the triangle inequality and Remark \ref{rk241}, we have
\begin{align}\label{ineq8}
& \left\Vert \int_s^t f(u)du- \int_0^{T} \mathds{1}_{[s,t]}f(u)du\right\Vert_E  \\ &\leq \left\Vert \tilde{m}_{s,t}(\mathscr{I}_f)-\tilde{m}_{s,t}(\mathbb{T}_n)\right\Vert_E+\left\Vert \tilde{m}_{0,T}(\tilde{\mathbb{T}}_n)-\tilde{m}_{0,T}(\mathscr{I}_{\mathds{1}_{[s,t]}f})\right\Vert_E\nonumber\\
&\leq (t-s)^{\frac{1}{2}}\left\Vert  \mathscr{I}_f- \mathbb{T}_n\right\Vert_{\gamma(0,T;E)}+ T^{\frac{1}{2}}\left\Vert  \tilde{\mathbb{T}}_n-\mathscr{I}_{\mathds{1}_{[s,t]}f}\right\Vert_{\gamma(0,T;E)}\nonumber\\
&\leq  2\, T^{\frac{1}{2}}\left\Vert  \mathscr{I}_f- \mathbb{T}_n\right\Vert_{\gamma(0,T;E)}.\nonumber 
\end{align}
Thus, by tending $n$ to infinity in (\ref{ineq8}) we get the remaining aim.
\end{proof}
The integrals of functions, as described in Proposition \ref{prop245}, can be compared with the Bochner integrals of these functions when the underlying space possesses one of the two geometric properties type $2$, cotype $2$ defined in [\cite{neerven2}, Definition $7.1.1$].\\
The reader can consult [\cite{BSEI}, Proposition $2.8$] for the proof of the following result.
%We may think of the space $L^2([s,t])\otimes E$, in another way as the space of square integrable functions defined on $[s,t]$ with finite-dimensional range in $E$, and for every $(h,e)\in L^2([s,t])\times E$, the Pettis integral operator of $u\longmapsto h(u)e$ will be equal to $h\otimes e$.\\ As a space of functions, $L^2([s,t])\otimes E$ is dense in $L^2([s,t];E)$ according to [\cite{neerven1}, Lemma $1.2.19$].
\begin{prop}\label{prop246}
Let $f:[0,T]\longrightarrow E$ be a $\lambda$-version of a measurable function.
\\If one of the two following assertions  holds:\vspace{0.2cm}\\
\hspace*{1cm} $(i)$\;\, $E$ has cotype $2$ and $f\in \gamma(0,T;E)$,\\
\hspace*{1cm} $(ii)$ $E$ has type $2$ and $f\in L^2([0,T];E)$,\vspace{0.2cm}\\	
then we have $$\int_s^tf(u)du=\int_{[s,t]}^{Bochner}f(u)du.$$
\end{prop}
\subsection{Stochastic Integration in UMD Spaces}
The class of UMD spaces was introduced in the 1970s by Maurey and Burkholder
and has since been studied by many authors. Examples of UMD spaces include Hilbert spaces, Lebesgue spaces $L^p(\mu)$ (for $1<p< \infty$), and certain Sobolev spaces. For further details, we refer the reader to \cite{UMD}.\\ Note that every UMD space is reflexive [\cite{neerven1}, Theorem 4.3.3].
\begin{defi}
A Banach space $F$ is said to be a UMD space if for all real numbers $q>1,$ there is a finite constant $\eta_{q,F}\geq 0$ (depending on $q$ and $F$) such that for all $F$-valued $L^{q}$-martingales $(\mathscr{M}_j)_{1\leq j\leq n}$ on a probability space $(\Omega',\mathscr{F}',\mathbb{P}')$, and for all choice of signs $\epsilon_j\in\left\lbrace-1,+1\right\rbrace,\,1\leq j\leq n,$ one has:
\begin{eqnarray}
\mathbb{E}'\left(\left\Vert \sum_{j=1}^n \epsilon_j (\mathscr{M}_j-\mathscr{M}_{j-1})\right\Vert_F^q\right) \leq \eta_{q,F}^q\mathbb{E}'\left(\left\Vert \sum_{j=1}^n (\mathscr{M}_j-\mathscr{M}_{j-1})\right\Vert_F^q\right),
\end{eqnarray}
where $\mathscr{M}_0=0_E$ by convention.
\end{defi}

In this subsection $E$ will be a separable UMD space. Now, let us suppose that our filtration
$\mathbb{F}:=\Big(\mathscr{F}_t\big)_{0\leq t\leq T}$ is generated by the Brownian motion $(W_t)_{t\in [0,T]}$, augmented by the Null-sets.\bigskip\\
We next define the stochastic integral of $E$-valued $\mathbb{F}$-adapted step processes:

\begin{defi}\label{def251}
	\begin{enumerate}
		\item[•]
A mapping $\varphi:[0,T]\times \Omega\longrightarrow E$ is called $\mathbb{F}$-adapted step process if it is of the following form
\begin{eqnarray*}
\varphi(t,\omega)=\mathds{1}_{\left\{0\right\}\times A}(t,\omega)\;\tilde{e}+\sum_{i=1}^n\sum_{j=1}^m\mathds{1}_{(t_{i-1},t_i]\times D_{ji}}(t,\omega)\;e_{ji},
\end{eqnarray*}
where $n,m$ are natural numbers; $0\leq t_0<...<t_n\leq T$; the set $A$ is from $\mathscr{F}_0$; the sets $D_{1i},...,D_{mi}$ are disjoints from $\mathscr{F}_{t_{i-1}}$, and the vectors $\tilde{e},e_{ji}$ are in $E$.
       \item[•]
The stochastic integral of the $\mathbb{F}$-adapted step process $\varphi$ of the previous form, with respect to $W$, is defined for each $s\in [0,T]$ as follows
$$\int_0^s \varphi_u dW_u=\sum_{i=1}^n\sum_{j=1}^m\mathds{1}_{D_{ji}}\left(W_{s\wedge t_i}-W_{s\wedge t_{i-1}}\right)\;e_{ji}.$$
	\end{enumerate}
\end{defi}
Now, we introduce the notion of $L^p$-stochastic  integrability.
\begin{defi}
A measurable and adapted process $\varphi:[0,T]\times \Omega\longrightarrow E$ is said to be $L^p$-stochastically integrable with respect to $W$ if the following conditions are fulfilled:
\begin{itemize}
\item[1)]For all $x^*\in E^*$ we have $\langle \varphi,x^*\rangle\in L^p(\Omega;L^2([0,T]))$.
\item[2)]There exists a sequence of\; $\mathbb{F}$-adapted step processes $\varphi^n:[0,T]\times \Omega\longrightarrow E$ such that:
\begin{itemize}
\item[i)]$(\varphi^n)_{n\geq 1}$ converges in measure on $[0,T]\times\Omega$ to $\varphi$,
\item[ii)]there exists a random variable $\mathbb{I}\in L^p(\Omega;E)$ such that
$$\lim_{n\rightarrow \infty}\int_0^T\varphi_u^ndW_u=\mathbb{I}\textrm{ in }L^p(\Omega;E).$$
\end{itemize}
\end{itemize}
In this case, $\mathbb{I}$ is called the stochastic integral of $\varphi$ with respect to $W$, and will be denoted by $\mathbb{I}:=\int_0^T\varphi_udW_u$.
\end{defi}
\begin{rk}
 It should be noted that, if the process $\varphi$ is $L^p$-stochastically integrable, then its restriction on $[0,s]\times \Omega$ is also $L^p$-stochastically integrable for every $s\in (0,T)$.\\
\end{rk}
Next, to invest an important characterization of the stochastic integrability, we need to explore the extension of the notion announced in Definition \ref{def241}:
\begin{defi}\label{def242}
	Let $\phi:[0,T]\times \Omega\longrightarrow E$ be a $\mathscr{B}([0,T])\otimes \mathscr{F}$-measurable process. $\phi$ is said to define a random variable $X:\Omega\longrightarrow \gamma(0,T;E)$ if the following conditions hold:
	\begin{eqnarray*}
		\left\{
		\begin{array}{lll}
			(i) & \forall x^*\in E^*,\;\langle \phi,x^*\rangle \in L^2([0,T])\textrm{ for a.s. }\, \omega\in\Omega
			;
			
			\\ (ii) & \forall h\in L^2([0,T]),\; \forall x^*\in E^*,\,\big \langle  X(\omega)h,x^*\big \rangle=\displaystyle\int_{[0,T]} h(u)\langle\phi(u,\omega),x^*\rangle d u\quad  \textrm{ for a.s. }\, \omega\in\Omega.
			
		\end{array}
		\right.
	\end{eqnarray*}
\end{defi}
\noindent We will write $\phi \in L^p\left(\Omega;\gamma(0,T;E)\right)$ instead of $X\in L^p\left(\Omega;\gamma(0,T;E)\right)$ and the same for its norm, where $\phi$ is a measurable process defining the random variable $X$.\\

%\begin{defi}
%	The Banach space $L_{\mathbb{F}}^p\left(\Omega;\gamma(0,T;E)\right)$ is defined as the closure in $L^p\left(\Omega;\gamma(0,T;E)\right)$ of the $E$-valued $\mathbb{F}$-adapted step processes.
%\end{defi}
\noindent Now, thanks to Theorem $3.6$ in \cite{neerven3}, a measurable and adapted $E$-valued mapping $\varphi:[0,T]\times \Omega\rightarrow E$ is $L^p$-stochastically integrable if and only if it defines an element of $L^p\left(\Omega;\gamma(0,T;E)\right)$. Moreover, we can confirm the existence of a sequence $(\varphi^n)$ of $E$-valued $\mathbb{F}$-adapted step processes converging to $\varphi$ in $L^p\left(\Omega;\gamma(0,T;E)\right)$, according to [\cite{neerven3}, Proposition $2.11$, Proposition $2.12$]. Additionally, we present the following result.
\begin{thm}[It\^o isomorphism]\label{thm251}
For any measurable and adapted process $\varphi\in L^p\left(\Omega;\gamma(0,T;E)\right)$, we have the following estimates
$$\mathbb{E}\left(\sup_{t\in [0,T]}\left\Vert \int_0^t\varphi_u dW_u\right\Vert^p\right)\eqsim_p\mathbb{E}\left(\left\Vert \int_0^T\varphi_u dW_u\right\Vert^p\right)\eqsim_{p,E}\left\Vert \varphi\right\Vert_{L^p\left(\Omega;\gamma(0,T;E)\right)}^p.$$
\end{thm}

\begin{rk} It should be noted that the UMD property for the space $E$ is necessary in Theorem \ref{thm251}. For further details, the reader is referred to \cite{neerven3}.
	\end{rk} 
Next, we analyze the behavior of the quantities generated by applying the integral studied in the previous subsection to a process defining an element of $L^p\left(\Omega;\gamma(0,T;E)\right)$.
\begin{prop}\label{prop247}
	Let $f:[0,T]\times \Omega\longrightarrow E$ be a measurable and adapted mapping defining an element of $L^p(\Omega;\gamma(0,T;E))$. Let $s,t$ be two real numbers from $[0,T]$ such that $s<t$.
	Then, the function $\omega\longmapsto \int_s^tf(u,\omega)du$ defines an element of $L^p(\Omega,\mathscr{F}_t;E)$ and we have
	$$\left\Vert \int_s^tf(u,.)du\right\Vert_{L^p(\Omega;E)}\leq (t-s)^{\frac{1}{2}}\left\Vert f\right\Vert_{L^p(\Omega;\gamma(0,T;E))}.$$
\end{prop}
\begin{proof}
	There exists a sequence $(f^n)_{n\geq 1}$ of $\mathbb{F}$-adapted step processes converging in $L^p(\Omega;\gamma(0,T;E))$ to $f$. Moreover, Lemma $2.7$ in \cite{neerven3} asserts that there exists a negligible set $\mathscr{N}$ from $\mathscr{F}_0$ such that for each $\omega$ in ($\Omega \setminus \mathscr{N}$), the mapping $f(., \omega)$ defines an element of $\gamma(0, T; E)$, and we have
	\begin{eqnarray}\label{ineq3}
		\left\Vert \int_s^tf(u,\omega)du-\int_s^tf^{\phi(n)}(u,\omega)du\right\Vert_E\leq (t-s)^{\frac{1}{2}}\left\Vert f(u,\omega)-f^{\phi(n)}(u,\omega)\right\Vert_{\gamma(0,T;E)},
	\end{eqnarray}
	where the subsequence $(f^{\phi(n)}(.,\omega))_{n\geq 1}$ converges to $f(.,\omega)$ in $\gamma(0,T;E)$ for every $\omega$ in ($\Omega \setminus \mathscr{N}$).\\
	Obviously, $\int_s^tf^{\phi(n)}(u,.)du$ is an $\mathscr{F}_t$-measurable simple function for every natural number $n$, thus by using (\ref{ineq3}) we deduce that $\int_s^tf(u,.)du$ is a $\mathbb{P}$-version of an $\mathscr{F}_t$-measurable random variable. Further, by the dominated convergence theorem, we can infer the remaining required result.
\end{proof}
\subsection{Set-valued Conditional Expectation, Set-valued Martingale}
The content of this subsection is related to multivalued analysis. For an overview, the reader may consult \cite{papa,Hiai,Him,set-stochastic-ch}.

Concerning the measurability of set-valued functions, there are several definitions in various situations. We deal with the following.\\
Let $(\mathcal{S},\mathscr{A},\mu)$ be a finite measure space, and $q\in [1,+\infty)$.
\begin{defi}
	Let $H:\mathcal{S}\longrightarrow \mathscr{K}\left( E\right)$ be a set-valued function.
	\begin{enumerate}
		\item[•]
		$H$ is called $\mathscr{A}$-measurable if for any open set $\mathscr{O}$ of $E$ we have
		\begin{eqnarray*}
			\left\lbrace s\in \mathcal{S} :\,H(s)\cap \mathscr{O}\neq \emptyset\right\rbrace \in \mathscr{A}.
		\end{eqnarray*}
		\item[•]
		For $\mathcal{S}=[0,T]\times \Omega$,
		\begin{itemize}
			\item[i)]$H$ is said to be a set-valued stochastic process if $H_t$ is a set-valued random variable for each $t\in [0,T]$, that is $H_t$ is $\mathscr{F}$-measurable for each $t\in [0,T]$. 
			\item[ii)]$H$ is called measurable if $H$ is $\mathscr{B}([0,T])\otimes \mathscr{F}$-measurable.
			\item[iii)]$H$ is said to be $\mathbb{F}$-adapted if $H_t$ is $\mathscr{F}_t$-measurable for every $t\in [0,T]$.
		\end{itemize}
	\end{enumerate}
\end{defi}
The notion of subtrajectory integrals plays an important role in the construction of many concepts related to multivalued analysis.
\begin{defi}\label{def}
	Let $G:\mathcal{S}\longrightarrow \mathscr{K}(E)$ be an $\mathscr{A}$-measurable set-valued function.
	\begin{itemize}
		\item[•]The subtrajectory integrals of $G$, $S_{\mathscr{A}}^q\left(G\right)$, is defined by  
%		and consists of all $g\in L^q\left(\Omega,\mathscr{A};E\right)$ such that $g(s)\in G(s)$ for almost everywhere $s\in \mathcal{S}$.
		\begin{eqnarray*}\label{sigma}
		S_{\mathscr{A}}^q\left(G\right)= \Big\{  g\in L^q\left(\mathcal{S},\mathscr{A};E\right) :g(s)\in G(s)\,\,\text{for almost everywhere} \,\,s\in \mathcal{S} \Big\}.
		\end{eqnarray*}
		\item[•]$G$ is said to be $q$-integrably bounded if there exists $\upsilon$ in $L^q\left(\mathcal{S},\mathscr{A}\right)$ such that for almost everywhere $s$ we have $$\left\vert G(s)\right\vert_{\varrho}\leq \upsilon(s).$$
	\end{itemize}
\end{defi}
We also introduce the concept of decomposability, which contributes to the development of various concepts in multivalued analysis, and is involved in some results presented in this paper.
\begin{defi}\label{def2}
	Let $\Delta$ be a nonempty subset of $L^q\left(\mathcal{S},\mathscr{A};E\right)$.
	\begin{itemize}
		\item[•]$\Delta$ is said to be decomposable with respect to $\mathscr{A}$, if for any finite $\mathscr{A}$-measurable partition $A_1,...,A_m$ of $\mathcal{S}$ and for any $h_1,...,h_m$ from $\Delta$, we have $\sum_{k=1}^m\mathds{1}_{A_k}h_k$ belongs to $\Delta$.
		\item[•]The decomposition hull $dec_{\mathscr{A}}(\Delta)$ $\big($resp. the closed decomposition hull $\overline{dec}_{\mathscr{A}}(\Delta)\big)$ of $\Delta$ is the smallest decomposable (resp. closed decomposable) subset of $L^q\left(\mathcal{S},\mathscr{A};E\right)$ containing $\Delta$.\\
		%We adopt $dec(\Delta),\,\overline{dec}(\Delta)$ instead of $dec_{\mathscr{A}}(\Delta),\,\overline{dec}_{\mathscr{A}}(\Delta)$ respectively, in no ambiguity situation.
	\end{itemize}
\end{defi}
Next, Theorem $3.8$ in [\cite{papa}, Chapter 2] and the decomposability of the subtrajectory integrals assert the existence of the following concept.
\begin{defi}
	Let $\mathscr{A}'$ be a sub-$\sigma$-algebra of $\mathscr{A}$ and $G:\mathcal{S}\longrightarrow \mathscr{K}(E)$ be an $\mathscr{A}$-measurable $q$-integrably bounded set-valued function. The set-valued conditional expectation of $G$ relative to $\mathscr{A}'$ is the $\mathscr{A}'$-measurable set-valued function from $\mathcal{S}$ into $\mathscr{K}(E)$, denoted by $\mathbb{E}^{\mathscr{A}'}\left(G\right)$, and satisfying $$S^q_{\mathscr{A}'}\left(\mathbb{E}^{\mathscr{A}'}\left(G\right)\right)=cl_{L^q\left(\mathcal{S},\mathscr{A}';E\right)}\left\lbrace \mathbb{E}^{\mathscr{A}'}\left(g\right):g\in S^q_{\mathscr{A}}\left(G\right)\right\rbrace,$$
	where $cl_{L^q\left(\mathcal{S},\mathscr{A}';E\right)}$ is the closure with respect to the strong topology of the space $L^q\left(\mathcal{S},\mathscr{A}';E\right)$.
\end{defi}
Now, we have sufficient ingredients to introduce the concept of a set-valued martingale.
\begin{defi} Let $\Psi:[0,T]\times \Omega\longrightarrow \mathscr{K}_c(E)$ be an $\mathbb{F}$-adapted set-valued function.\\ Assume that $\Psi_T$ is $q$-integrably bounded.
%	\begin{eqnarray*}
%	\exists\upsilon\in L^p(\Omega,\mathscr{F}_T),\quad\sup_{t\in [0,T]} \left\vert \Psi_t\right\vert_{\varrho}\leq \upsilon\quad a.s.\;\omega\in\Omega.	
%	\end{eqnarray*}
	\begin{itemize}
		\item $\Psi$ is said to be a set-valued martingale if for every $0\leq t_1\leq t_2\leq T$, $\mathbb{E}^{\mathscr{F}_{t_1}}\left(\Psi_{t_2} \right)=\Psi_{t_1}\;a.s.\;\omega$.
		\item An $E$-valued $\mathbb{F}$-martingale $(\psi)_{0\leq t\leq T}$ is called a martingale selection of $\Psi$ if $\psi_t$ belongs to $S^p_{\mathscr{F}_t}\left(\Psi_t\right)$ for every $t$ in $[0,T]$.
	\end{itemize}
\end{defi}
\subsection{Hukuhara Difference}
In the metric space $(\mathscr{K}_{cwcmpt}(E), \varrho)$, the Minkowski sum and scalar multiplication of sets are two operations. However, they generally do not allow for the construction of an opposite element.\\Hukuhara, in \cite{hukuhara}, introduced a notion of set difference for certain nonempty compact convex sets in  $\mathbb{R}^n$, which has been used to study specific differential inclusions, as discussed in [\cite{set-stochastic-ch}, Chapter 7]. 
Based on this, we extend the concept as follows.
\begin{defi}\label{defi201}
Let $A,B$ be two elements from $\mathscr{K}_{cwcmpt}(E)$. 
The nonempty convex closed subset $C$ of $E$ is said to be the Hukuhara difference of $A$ and $B$ if it satisfies $B + C = A$. In case where $C$ exists, it will be denoted by $A\circleddash B$.
\end{defi}
\begin{rk}
We notice that, in case of existence, the Hukuhara difference $A\circleddash B$ is unique by virtue of the cancellation law lemma [\cite{papa},  Chapter 1, Lemma $2.74$], and it is easily to check that
\begin{eqnarray}\label{formula1}
\left\vert A\circleddash B\right\vert_{\varrho}\leq \left\vert A\right\vert_{\varrho}+\left\vert B\right\vert_{\varrho},
\end{eqnarray}
where  $\left\vert\, .\,\right\vert_{\varrho}$ is defined by (\ref{n.hausdorff}).
\end{rk}
Building on Proposition 4.16 in \cite{Hamel}, we present a key result that facilitates the proof of the measurability of the Hukuhara difference between two set-valued functions.
\begin{prop}\label{prop203}
Assume that $E$ is a reflexive space.
Let $A$ and $B$ be two nonempty convex weakly compact subsets of $E$, and let $(x^*_n)_{n\geq 1}$ be a sequence in $E^*\setminus\left\lbrace 0_{E^*}\right\rbrace$ that is dense in $E^*$.\\If $A\circleddash B$ exists, then we have
$$A\circleddash B=\bigcap_{n\geq 1}\left\lbrace x\in E:\langle x,x^*_n\rangle\,\,\leq\,\, \sup_{a\in A}\langle a,x^*_n\rangle-\sup_{b\in B}\langle b,x^*_n\rangle\right\rbrace.$$
\end{prop}
\begin{proof}
Let $x\in A\circleddash B$. Since $x+B\subset A$, then one can get obviously
\begin{eqnarray}\label{assertion4}
\forall x^*\in E^*\setminus\left\lbrace 0_{E^*}\right\rbrace,\; \langle x,-x^*\rangle\,\,\geq\,\, \inf_{a\in A}\langle a,-x^*\rangle-\inf_{b\in B}\langle b,-x^*\rangle.
\end{eqnarray}
Conversely, take $x\in E$ such that assertion (\ref{assertion4}) holds, and assume that $x\notin A\circleddash B$.\\
Let us choose $b\in B$ such that $x+b\notin A$. Since $A$ is nonempty closed convex subset of the Banach space $E$, then by applying hyperplane separation theorem we derive the existence of $x^*\in E^*\setminus\left\lbrace 0_{E^*}\right\rbrace$ satisfying $\langle x+b,-x^*\rangle\;<\inf_{a\in A}\langle a,-x^*\rangle$, which is contradictory with (\ref{assertion4}). Therefore,
$$A\circleddash B=\bigcap_{x^*\in E^*\setminus\left\lbrace 0_{E^*}\right\rbrace}\left\lbrace x\in E:\langle x,x^*\rangle\leq \sup_{a\in A}\langle a,x^*\rangle-\sup_{b\in B}\langle b,x^*\rangle\right\rbrace.$$
It remains to use the continuity of the support functions $\sup_{a\in A}\langle a,.\rangle$ and $\sup_{b\in B}\langle b,.\rangle$ to conclude the result.
\end{proof} 
%%%%%%%%%%%%%%%%%%%%%%%%%%%%%%%%%%%%%Section4%%%%%%%%%%%%%%%%%%%%%%%%%%%%%%%%%%%%%%%%%%%%%%%

\section{The Revised Set-valued Stochastic Integrals in UMD Spaces}
Let $E$ be a separable UMD space, and assume that $\mathbb{F}$ is the natural filtration generated by a standard Brownian motion $W=\left(W_t\right)_{t\in [0,T]}$ augmented by all $\mathbb{P}$-null sets of $\mathscr{F}$. Let us denote by 
$\mathcal{L}_{\mathbb{F}}^p([0,T]\times \Omega;E)$ the set of (equivalence classes of) measurable and adapted mappings from $[0,T]\times \Omega$ into $E$ which are $L^p$-stochastically integrable.\bigskip\\
Here, it should be noted that two measurable processes $\phi_1,\phi_2:[s,t]\times\Omega\rightarrow E$, define the same random variable $X:\Omega\rightarrow\gamma(s,t;E)$ if and only if $\phi_1(u,\omega)=\phi_2(u,\omega)$ for almost all $(u,\omega)\in [s,t]\times\Omega$. Conversely, two random variables $X_1, X_2 :\Omega \rightarrow\gamma(s, t; E)$ are defined by the same measurable process $\phi$ if and only if $X_1(\omega) = X_2(\omega)$ for a.s. $\omega \in \Omega$.\bigskip\\
Let $L_0^p(\Omega,\mathscr{F}_T;E)$ be the closed subspace of $L^p(\Omega,\mathscr{F}_T;E)$ consisting of all $\mathscr{F}_T$-measurable random variables with mean zero.\\
We begin by determining the essential topological structure of the space $\mathcal{L}_{\mathbb{F}}^p([0,T]\times \Omega;E)$.
\begin{prop}\label{prop301}
The space $\mathcal{L}_{\mathbb{F}}^p([0,T]\times \Omega;E)$ endowed with the norm $$f\in \mathcal{L}_{\mathbb{F}}^p([0,T]\times \Omega;E)\longmapsto\left\Vert f\right\Vert_{L^p(\Omega;\gamma(0,T;E))},$$ is a Banach space.
\end{prop}
\begin{proof}
Let $(f^n)_{n\geq 1}$ be a Cauchy sequence from $\mathcal{L}_{\mathbb{F}}^p([0,T]\times \Omega;E)$. According to Theorem \ref{thm251}, we have, for every $m, n \in \mathbb{N}$, that
\begin{eqnarray*}
\left\Vert f^m-f^n\right\Vert_{L^p(\Omega;\gamma(0,T;E))}\eqsim_{p,E} \left\Vert \int_0^T \left(f_u^m-f_u^n\right) dW_u\right\Vert_{L_0^p(\Omega,\mathscr{F}_T;E)}.
\end{eqnarray*}
Thus, $\left(\int_0^T f_u^n dW_u\right)_{n\geq 1}$ converges to a random variable $h$ in $L_0^p(\Omega,\mathscr{F}_T;E)$. According to the representation of Brownian $L^p$-martingales [\cite{BSEE.B}, $(3.2)$], we conclude that there exists $f$ belonging to $\mathcal{L}_{\mathbb{F}}^p([0,T]\times \Omega;E)$ such that
\begin{eqnarray*}
h=\int_0^T f_u dW_u\quad\text{in}\; L^p(\Omega;E).
\end{eqnarray*}
Therefore, $(f_n)_{n\geq 1}$ converges to $f$ in $\mathcal{L}_{\mathbb{F}}^p([0,T]\times \Omega;E)$. This completes the proof.
\end{proof}
In order to define our revised set-valued stochastic integral, we first define and study an associated concept called the set-valued functional stochastic integrals with increments.
\subsection{Set-valued Functional Stochastic Integrals with Increments}

Let $\mathfrak{S}$ be a nonempty bounded convex subset of $E \times \mathcal{L}_{\mathbb{F}}^p([0,T] \times \Omega; E)$, and let $t \in (0, T]$.\\ Let $cl(\mathfrak{S})$ denote the closure in the space $E \times \mathcal{L}_{\mathbb{F}}^p([0, T] \times \Omega; E)$ endowed with the norm $$(x,f)\longmapsto  \left\Vert x\right\Vert_E+\left\Vert f\right\Vert_{L^p(\Omega;\gamma(0,T;E))}.$$
% defined by:
%$$(x,f) \longmapsto  \left\Vert x\right\Vert_E+\left\Vert f\right\Vert_{L^p(\Omega;\gamma(0,T;E))} .$$
\begin{defi}
The set-valued functional stochastic integral of $\mathfrak{S}$ over $[0,t]$ is defined as follows:
\begin{eqnarray}
\mathbb{I}_t(\mathfrak{S}):=\left\lbrace x+\int_0^tf_udW_u:(x,f)\in \mathfrak{S}\right\rbrace.
\end{eqnarray}
\end{defi}
\noindent Let $proj_1$ be the projection defined on $E\times \mathcal{L}_{\mathbb{F}}^p([0,T]\times \Omega;E)$ onto the first coordinate set. \\

 {The following theorem  is a generalization of the first property of Lemma $3.1$ in \cite{kis2013}.}
\begin{thm}\label{thm303}The set $\mathbb{I}_t(\mathfrak{S})$ is a bounded subset of $L^p(\Omega,\mathscr{F}_t;E)$, and we have
	\begin{eqnarray}\label{eq301}
	cl_{L^p}\left(\mathbb{I}_t(\mathfrak{S})\right)=\mathbb{I}_t\left(cl(\mathfrak{S})\right),	
	\end{eqnarray} 
	where "$cl_{L^p}$" denotes the closure with respect to the norm topology of the space $L^p\left(\Omega,\mathscr{F}_t;E\right)$.
\end{thm}
\begin{proof}
	\begin{itemize}
		\item [$\bullet$]Let $\alpha$ and $\beta$ two positive real numbers such that $\alpha +\beta =1$. We have the following
		\begin{eqnarray*}
			\alpha \mathbb{I}_t(\mathfrak{S})+\beta \mathbb{I}_t(\mathfrak{S})&=&\left\lbrace \alpha x+\beta y+\int_0^t(\alpha f_u+\beta g_u)dW_u:(x,f),(y,g)\in\mathfrak{S}\right\rbrace\\
			&=&\left\lbrace z+\int_0^t h_udW_u:(z,h)\in(\alpha\mathfrak{S}+\beta \mathfrak{S})\right\rbrace.
		\end{eqnarray*}
		Hence, the convexity of $\mathbb{I}_t(\mathfrak{S})$ is derived from that of $\mathfrak{S}$.
		\item [$\bullet$]Let $\left((x^n, f^n)\right)_{n \geq 1}$ be a sequence in $\mathfrak{S}$ such that $\left(\mathbb{I}_t(x^n, f^n)\right)_{n \geq 1}$ converges in $L^p(\Omega, \mathscr{F}_t; E)$.
		By applying expectation, we conclude that $(x^n)_{n \geq 1}$ is a Cauchy sequence in $E$, thus it converges in $E$, say to $x$. On the other hand, the mapping:
		\begin{eqnarray}
			\begin{array}{ccccc}
				 &  & E\times \mathcal{L}_{\mathbb{F}}^p([0,T]\times \Omega;E) & \longrightarrow & L^p\left(\Omega,\mathscr{F}_T;E\right)\\
				& & (z,g) & \longmapsto & z+\displaystyle\int_0^Tg_udW_u,\\
			\end{array}\label{isomor2}
		\end{eqnarray}
		 establishes an isomorphism from $E \times \mathcal{L}^p_{\mathbb{F}}\left([0, T] \times \Omega; E\right)$ onto $L^p(\Omega, \mathscr{F}_T; E)$, and $L^p(\Omega, \mathscr{F}_T; E)$ is indeed a reflexive Banach space. Thus, $E \times \mathcal{L}^p_{\mathbb{F}}\left([0, T] \times \Omega; E\right)$ is itself reflexive by virtue of Proposition 1.11.8 in \cite{reflexive}. Then, by virtue of Lemma 1.13.3 in \cite{reflexive}, we have a weakly convergent subsequence $((x^{\varphi(n)}, f^{\varphi(n)}))_{n \geq 1}$ converging to an element $(y, f)$ in $E \times \mathcal{L}^p_{\mathbb{F}}\left([0, T] \times \Omega; E\right)$. Additionally, $(y, f) \in cl(\mathfrak{S})$ by virtue of to Theorem 10.1.2 in \cite{topog}.
		 
		 Since $proj_1$ is a continuous linear mapping,  it follows from Theorem 10.2.2 in \cite{topog} that the subsequence $(x^{\varphi(n)})_{n \geq 1}$ weakly converges to $y$. Henceforth, $y = x$.
		 \\Thanks to It\^o isomorphism theorem and Remark \ref{rk241}, we have
		 \begin{eqnarray}
		 \left\Vert \mathbb{I}_t(z, g)\right\Vert_{L^p(\Omega,\mathscr{F}_t;E)}&\leq& \left\Vert z\right\Vert_E+	\left\Vert \int_0^tg_udW_u\right\Vert_{L^p(\Omega,\mathscr{F}_T;E)}\nonumber \\ &\quad\,\,\, \eqsim_{p,E}&\left\Vert z\right\Vert_E+\left\Vert \mathds{1}_{[0,t]}g\right\Vert_{L^p(\Omega;\gamma(0,T;E))}\label{ineq7}\\
		 &\leq &\left\Vert z\right\Vert_E+\left\Vert g\right\Vert_{L^p(\Omega;\gamma(0,T;E))},\nonumber 
		 \end{eqnarray} 
		  for any $(z,g)\in E\times \mathcal{L}^p_{\mathbb{F}}\left([0, T] \times \Omega; E\right)$. It turns out that $(\int_0^tf_u^{\varphi(n)}dW_u)_{n\geq 1}$ is weakly convergent to $\int_0^tf_udW_u$.
		Therefore, $\lim_{n\rightarrow \infty}\mathbb{I}_t(x^n,f^{n})=\mathbb{I}_t(x,f)$ in $L^p(\Omega,\mathscr{F}_t;E)$, which asserts that $cl_{L^p}\left(\mathbb{I}_t(\mathfrak{S})\right)$ is contained in $\mathbb{I}_t\left(cl(\mathfrak{S})\right)$. By continuity, we derive the reverse inclusion, and thus, equality (\ref{eq301}) follows.
		\item [$\bullet$] Using (\ref{ineq7}), we conclude that the boundedness of $\mathfrak{S}$ in $E\times\mathcal{L}_{\mathbb{F}}^p([0,T]\times \Omega;E)$ implies the boundedness of $\mathbb{I}_t\left(\mathfrak{S}\right)$ in $L^p(\Omega,\mathscr{F}_t;E)$.
	\end{itemize}
Finally, the proof is accomplished.
\end{proof}
\begin{rk}
 From (\ref{isomor2}), we note that a set belongs to $\mathscr{K}_{cwcmpt}\left(E \times \mathcal{L}_{\mathbb{F}}^p([0,T] \times \Omega; E)\right)$ if and only if it is a nonempty, bounded, closed, and convex subset of the space $E \times \mathcal{L}_{\mathbb{F}}^p([0,T] \times \Omega; E)$.
\end{rk}
%\begin{rk}
%If $\mathfrak{S}$ is closed, then $\overline{dec}\left\lbrace \mathbb{I}_t(\mathfrak{S})\right\rbrace$ is convex and weakly closed in $L^p(\Omega,\mathscr{F}_t;E)$. Indeed, due to Theorem \ref{thm303} we have $\mathbb{I}_t(\mathfrak{S})$ is a closed convex subset of $L^p(\Omega,\mathscr{F}_t;E)$, and by taking into account Lemma $3.3.3$ in \cite{set-stochastic-ch} we infer that $\overline{dec}\left(\mathbb{I}_t(\mathfrak{S})\right)$ is a closed convex subset of $L^p(\Omega,\mathscr{F}_t;E)$, thus the aim follows.
%\item[•]For $t=T$ and $\mathfrak{S}=\left\lbrace 0_E\right\rbrace\times \mathscr{G}$, where $\mathscr{G}$ is a nonempty subset of $\mathcal{L}_{\mathbb{F}}^p([0,T]\times \Omega;E)$, the equivalence between the boundedness (resp. closedness) of $\mathscr{G}$ and the boundedness (resp. closedness) of $\mathbb{I}_T(\mathfrak{S})$ follows by virtue of It\^o isomorphism theorem.
%\end{rk}

The two remaining results require the separability of the probability space.
In the first result, we may provide a sequence that represents, using the weak closedness, both the set-valued functional stochastic integral with increments of the set $\mathfrak{S}$, and its closed decomposition hull, which is a member of the basic equality in Definition \ref{defi301}.
%\\For further clarity, we refer the reader to \cite{Brown}, page $174$, for a detailed explanation of separability. It is also worth mentioning that if the probability space is separable, its completion will be separable as well.
\begin{thm}\label{thm304}
If $\mathfrak{S}$ is closed, then $\mathbb{I}_t(\mathfrak{S})$ is a weakly compact subset of $L^p(\Omega,\mathscr{F}_t;E)$. Furthermore, assuming that  $\left(\Omega,\mathscr{F},\mathbb{P}\right)$ is separable, we can select a  sequence $\left((x^n,f^n)\right)_{n\geq 1}$ from $\mathfrak{S}$ such that
\begin{itemize}
\item[(i)]$\mathbb{I}_t(\mathfrak{S})=cl_w\left\lbrace x^n+\displaystyle\int_0^tf^n_u dW_u:n\geq 1\right\rbrace$,
\item[(ii)]$\overline{dec}_{\mathscr{F}_t}(\mathbb{I}_t(\mathfrak{S}))=cl_w\left(dec_{\mathscr{F}_t}\left\lbrace x^n+\displaystyle\int_0^tf^n_u dW_u:n\geq 1\right\rbrace\right)$,
\end{itemize}
where "$cl_w$" denotes the closure with respect to the weak topology of the space $L^p(\Omega,\mathscr{F}_t;E)$.
\end{thm}
\begin{proof}
\begin{itemize}
\item Theorem \ref{thm303} guaranties that $\mathbb{I}_t(\mathfrak{S})$ is a convex closed and bounded subset of $L^p(\Omega,\mathscr{F}_t;E)$. Thus, due to Theorem 1.2.3 in \cite{set-stochastic-ch} we derive the weak-compactness of $\mathbb{I}_t(\mathfrak{S})$ since $E$ is reflexive.

For the remainder of the proof, we assume that $\left(\Omega,\mathscr{F},\mathbb{P}\right)$ is separable.
\item  Since $L^p\left(\Omega,\mathscr{F}_t;E\right)$ is a separable reflexive Banach space, then its associated topological dual is also separable. Thus, the unit ball $\overline{\mathcal{B}}$ of $L^p\left(\Omega,\mathscr{F}_t;E\right)$ is metrizable with respect to the weak topology. Besides the weak-compactness property given by Banach-Alaoglu-Bourbaki theorem we conclude that $\overline{\mathcal{B}}$, endowed with the induced weak topology of $L^p\left(\Omega,\mathscr{F}_t;E\right)$, is a separable topological space.\\
On the other hand there is a positive real number $\beta$ such that $\beta\mathbb{I}_t(\mathfrak{S})$ be contained in $\overline{\mathcal{B}}$, thus by combining with the weak closedness of the convex set $\mathbb{I}_t(\mathfrak{S})$ we derive the existence of a sequence $\left((x^n,f^n)\right)_{n\geq 1}$ which satisfies the property $(i)$.
\item  Let $h$ be from $cl_w\left(dec_{\mathscr{F}_t}\left(cl_w\left\lbrace x^n+\int_0^tf_u^ndW_u:n\geq 1\right\rbrace\right)\right)$ where we still adopt the same sequence $\left((x^n,f^n)\right)_{n\geq 1}$ involved in the first property.
\\Since $h$  has a countable neighbourhood basis, then we may choose a sequence $(h^n)_{n\geq 1}$ in $L^p\left(\Omega,\mathscr{F}_t;E\right)$ weakly converging to $h$. Let for each natural number $n$, $(\Omega_{n,l})_{1\leq l\leq m_n}$ be an $\mathscr{F}_t$-measurable partition of $\Omega$ and $(h^{n,l})_{1\leq l\leq m_n}$ be a sequence from $cl_w\left\lbrace x^n+\int_0^tf_u^ndW_u:n\geq 1\right\rbrace$ with $h^n=\sum_{l=1}^{m_n}\mathds{1}_{\Omega_{n,l}}h^{n,l}$.
\\Let $\left((x^{\varphi_{n,l}(d)},f^{\varphi_{n,l}(d)})\right)_{d\geq 1}$ be a subsequence such that $\left(x^{\varphi_{n,l}(d)}+\int_0^tf_u^{\varphi_{n,l}(d)}dW_u\right)_{d\geq 1}$ is weakly convergent to $h^{n,l}$ for every pair $(n, l)$ of natural numbers with $1\leq l\leq m_n$.
\\From the boundedness of the sequence $\left(\int_0^T f_u^{\varphi_{n,l}(d)}dW_u\right)_{d\geq 1}$ in the reflexive space $L^p(\Omega,\mathscr{F}_T;E)$, we derive due to It\^o isomorphism theorem that there is a subsequence $\left(f^{\varphi'_{n,l}\circ\varphi_{n,l}(d)}\right)_{d\geq 1}$ weakly convergent in $\mathcal{L}_{\mathbb{F}}^p([0,T]\times \Omega;E)$ to $f^{\varphi'_{n,l}\circ\varphi_{n,l}}$.\\ Similarly, $\left(x^{\varphi_{n,l}(d)}\right)_{d\geq 1}$ is bounded in $E$, then we may take $\left(x^{\varphi'_{n,l}\circ\varphi_{n,l}(d)}\right)_{d\geq 1}$ as a subsequence weakly convergent in $E$ to $x^{\varphi'_{n,l}\circ\varphi_{n,l}}$.\bigskip\\
Let $g$ be in $L^q\left(\Omega,\mathscr{F}_t;E^*\right)$, where $q$ is the conjugate index of $p$. Due to the following inequality
\begin{align*}
&\left\vert \mathbb{E}\left(\langle h,g\rangle\right)-\mathbb{E}\left(\left\langle\sum_{l=1}^{m_n}\mathds{1}_{\Omega_{n,l}}\left(x^{\varphi'_{n,l}\circ\varphi_{n,l}}+\int_0^tf_u^{\varphi'_{n,l}\circ\varphi_{n,l}}dW_u\right),g\right\rangle\right)\right\vert\\
&\leq  \left\vert \mathbb{E}\left(\left\langle h-h^n,g\right\rangle\right)\right\vert+\left\vert \mathbb{E}\left(\sum_{l=1}^{m_n}\left\langle h^{n,l}-x^{\varphi'_{n,l}\circ\varphi_{n,l}(d)}-\int_0^tf_u^{\varphi'_{n,l}\circ\varphi_{n,l}(d)}dW_u,\mathds{1}_{\Omega_{n,l}}g\right\rangle\right)\right\vert\\
&+ \left\vert \mathbb{E}\left(\sum_{l=1}^{m_n}\left\langle x^{\varphi'_{n,l}\circ\varphi_{n,l}(d)}-x^{\varphi'_{n,l}\circ\varphi_{n,l}}+\int_0^tf_u^{\varphi'_{n,l}\circ\varphi_{n,l}(d)}dW_u-\int_0^t f_u^{\varphi'_{n,l}\circ\varphi_{n,l}}dW_u,\mathds{1}_{\Omega_{n,l}}g\right\rangle\right)\right\vert,
\end{align*}
we conclude after tending $d$ to infinity that for every $n$ we have
$$\left\vert \mathbb{E}\left(\langle h,g\rangle\right)-\mathbb{E}\left(\left\langle\sum_{l=1}^{m_n}\mathds{1}_{\Omega_{n,l}}\left(x^{\varphi'_{n,l}\circ\varphi_{n,l}}+\int_0^tf_u^{\varphi'_{n,l}\circ\varphi_{n,l}}dW_u\right),g\right\rangle\right)\right\vert\leq \left\vert \mathbb{E}\left(\langle h-h^n,g\rangle\right)\right\vert.$$
Thus, the sequence $\left(\sum_{l=1}^{m_n}\mathds{1}_{\Omega_{n,l}}\left(x^{\varphi'_{n,l}\circ\varphi_{n,l}}+\int_0^tf_u^{\varphi'_{n,l}\circ\varphi_{n,l}}dW_u\right)\right)_{n\geq 1}$ is weakly convergent to $h$ in the space $L^p\left(\Omega,\mathscr{F}_t;E\right)$. Therefore,
$$cl_w\left(dec_{\mathscr{F}_t}\left(cl_w\left\lbrace x^n+\int_0^tf_u^ndW_u:n\geq 1\right\rbrace\right)\right)=cl_w\left(dec_{\mathscr{F}_t}\left\lbrace x^n+\int_0^tf_u^ndW_u:n\geq 1\right\rbrace\right).$$
On the other hand, the convexity of the set $\mathbb{I}_t(\mathfrak{S})$ implies that $dec_{\mathscr{F}_t}(\mathbb{I}_t(\mathfrak{S}))$ is also convex. Therefore,
\begin{eqnarray*}
\overline{dec}_{\mathscr{F}_t}(\mathbb{I}_t(\mathfrak{S}))=cl_{L^p\left(\Omega,\mathscr{F}_t;E\right)}\left(dec_{\mathscr{F}_t}(\mathbb{I}_t(\mathfrak{S}))\right)=cl_w\left(dec_{\mathscr{F}_t}(\mathbb{I}_t(\mathfrak{S}))\right).
\end{eqnarray*}
\end{itemize}
Consequently, both goals have been achieved.
\end{proof}
%Now, we can stipulate that the  sequence $(x^n,f^n)$ can be chosen independent of $t$.
Now, by considering the closure with respect to the strong topology of the space  $L^p\left(\Omega,\mathscr{F}_t;E\right)$,  we can stipulate that the sequence $(x^n,f^n)$ may be chosen independently of $t$.

\begin{thm}\label{thm305}
If $(\Omega,\mathscr{F},\mathbb{P})$ is separable, then there exists $\left((x^n,f^n)\right)_{n\geq 1}\subset\mathfrak{S}$ such that for every $s\in (0,T]$ we have
$$cl_{L^p}\left(\mathbb{I}_s\left(\mathfrak{S}\right)\right)=cl_{L^p}\left\lbrace x^n+\int_0^sf^n_u dW_u:n\geq 1\right\rbrace,$$
where "$cl_{L^p}$" is the colsure with respect to the strong topology of the space  $L^p\left(\Omega,\mathscr{F}_s;E\right)$.
\end{thm}
\begin{proof}
Since $(\Omega,\mathscr{F},\mathbb{P})$ is separable and $\gamma(0,T;E)$ is a separable Banach space, it follows that $L^p(\Omega,\gamma(0,T;E))$ is also separable. Furthermore, $\mathcal{L}^p_{\mathbb{F}}\left([0,T]\times \Omega;E\right)$ is homeomorphic to a closed subset of $L^p(\Omega,\gamma(0,T;E))$; which indicates that it is a separable metric space. As a result, the Cartesian product  $E\times\mathcal{L}_{\mathbb{F}}^p([0,T]\times \Omega;E)$ is also separable. By inducing this topology on $\mathfrak{S}$ to $\mathcal{T}_{d_T}$, we get a separable subspace of $E\times\mathcal{L}^p_{\mathbb{F}}\left([0,T]\times \Omega;E\right)$. It follows that there exists a sequence $\left((x^n,f^n)\right)_{n\geq 1}$ such that $\mathfrak{S}=cl_{\mathcal{T}_{d_T}}\left\lbrace (x^n,f^n):n\geq 1\right\rbrace$.
Subsequently
\begin{eqnarray*}
cl(\mathfrak{S})=cl\left(\mathfrak{S}\cap cl\left\lbrace (x^n,f^n):n\geq 1\right\rbrace\right)\subset cl(\mathfrak{S})\cap cl\left\lbrace (x^n,f^n):n\geq 1\right\rbrace.
\end{eqnarray*}
Then, $cl(\mathfrak{S})=cl\left\lbrace (x^n,f^n):n\geq 1\right\rbrace$. According to Theorem \ref{thm303}, we conclude that
\begin{eqnarray*}
cl_{L^p}\left(\mathbb{I}_s(\mathfrak{S})\right)=\mathbb{I}_s\left(cl(\mathfrak{S})\right)=\mathbb{I}_s\left(cl\left\lbrace (x^n,f^n):n\geq 1\right\rbrace\right)=cl_{L^p}\left\lbrace x^n+\int_0^sf^n_u dW_u:n\geq 1\right\rbrace.
\end{eqnarray*}
Finally, the proof is accomplished.
\end{proof}
%\begin{rk}
%Let $t\in (0,T]$. For $\mathfrak{S}=\left\lbrace 0_E\right\rbrace\times \mathscr{G}$, where $\mathscr{G}$ is a nonempty subset of $\mathcal{L}_{\mathbb{F}}^p([0,T]\times \Omega;E)$, we have $\mathbb{I}_t(\mathfrak{S})$ is decomposable if and only if it is a singleton. Indeed, let $A$ belongs to $\mathscr{F}_t$, $f,g$ be from $\mathfrak{S}$, and assume $\mathbb{I}_t(\mathfrak{S})$ is decomposable. We have $\mathds{1}_A \int_0^tf_udW_u+\mathds{1}_{\Omega\setminus A}\int_0^tg_udW_u$ remains in $\mathbb{I}_t(\mathfrak{S})$. Thus, $\mathbb{E}\left(\mathds{1}_A\int_0^t(f_u-g_u)dW_u\right)=0$, and $\int_0^tf_udW_u=\int_0^t g_udW_u$ in $L^p(\Omega,\mathscr{F}_t;E)$.
%\end{rk}
\subsection{The Revised Set-valued Stochastic Integrals}
It should be noted that the existence of the following concept is attributed to Theorem $3.8$ in [\cite{papa}, Chapter $2$].
\begin{defi}\label{defi301} The $\mathscr{F}_t$-measurable set-valued $\Theta_t:\Omega\mapsto\mathscr{K}(E)$, which satisfies
 $$S^p_{\mathscr{F}_t}(\Theta_t)=\overline{dec}_{\mathscr{F}_t}\left(\mathbb{I}_t(\mathfrak{S})\right),$$ is called the revised set-valued stochastic integral of $\mathfrak{S}$ on $[0,t]$, and denoted by $\displaystyle\int_{[0,t]}^{\mathscr{R}}\mathfrak{S}\, dW_u$. 
\end{defi}

%
%\begin{rk}
% It should be noted that the revised set-valued stochastic integral is closed convex valued, and we have
%	\begin{eqnarray*}
%		\int_{[0,t]}^{\mathscr{R}}cl(\mathfrak{S})\, dW_u=\int_{[0,t]}^{\mathscr{R}}\mathfrak{S}\, dW_u\;a.s.\;\omega\in\Omega.	
%	\end{eqnarray*}
%%Additionally, the set-valued integrals defined by (\ref{g-s-v-int}) are also considered as a particular case.
%\end{rk}	

Theorem \ref{thm304} explicitly provides the subtrajectory integrals of the revised set-valued stochastic integral using a sequence contained in the integrated subset. Clearly, its existence depends on $t$.\\
Let $\mathscr{M}^{c,p}(E)$ denote the set of all $E$-valued $p$-martingales with respect to the filtration $\mathbb{F}$ that have continuous trajectories.
\\To provide a Casting representation of the studied integral, constructed using a sequence in $\mathfrak{S}$, we first state the following lemma, which is an important tool to prove that the set $$\left\lbrace\left(x+\int_0^tf_udW_u\right)_{T\geq t\geq 0},\;(x,f)\in \mathfrak{S}\right\rbrace$$ equipped with the topology induced by $\mathscr{M}^{c,p}(E)$ is separable.
 
\begin{lem}\label{lem301}
If the space $\left(\Omega,\mathscr{F},\mathbb{P}\right)$ is separable, then the real vector space $\mathscr{M}^{c,p}(E)$ is a separable closed subset of $L^p\left(\Omega;C([0,T],E)\right)$.
\end{lem}
\begin{proof}
Let $\left(M^n\right)_{n\geq 1}$ be a sequence from $\mathscr{M}^{c,p}(E)$ tending to $M$ in $L^p\left(\Omega;C([0,T],E)\right)$.\\Let $s,t$ be two real numbers with $0\leq s\leq t\leq T$. By the contractiveness of the conditional expectation we have for each natural number $n$
\begin{eqnarray*}
\left\Vert M_s-\mathbb{E}^{\mathscr{F}_s}(M_t)\right\Vert_{L^p(\Omega;E)}&\leq & \left\Vert M_s-M^n_s\right\Vert_{L^p(\Omega;E)}+\left\Vert \mathbb{E}^{\mathscr{F}_s}(M^n_t)-\mathbb{E}^{\mathscr{F}_s}(M_t)\right\Vert_{L^p(\Omega;E)}\\
&\leq & 2 \left\Vert M-M^n\right\Vert_{L^p\left(\Omega;C([0,T],E)\right)}.
\end{eqnarray*}
Thus, $M_s=\mathbb{E}^{\mathscr{F}_s}(M_t)$ in $L^p(\Omega;E)$, and $M$ is a $p$-integrable continuous martingale.
\\To prove the separability of $\mathscr{M}^{c,p}(E)$ we start by examining some particular subspaces.\\
Let $\left\Vert .\right\Vert_{\mathscr{M}^{c,p}(E)}$ be the induced norm on the vector subspace $\mathscr{M}^{c,p}(E)$ by that of $L^p\left(\Omega,\mathscr{F}_T;C([0,T],E)\right)$.\bigskip
\\We will first focus on the vector space $\mathscr{M}^{c,p}(\mathbb{R}) \otimes E$, consisting of all finite combinations of the form $\sum_{j=1}^m \phi^j x_j$, where $\phi^j \in \mathscr{M}^{c,p}(\mathbb{R})$ and $x_j \in E$.\\
Since $\left(\Omega,\mathscr{F},\mathbb{P}\right)$ is separable and $C([0,T],\mathbb{R})$ is a separable Banach space, then  $L^p\left(\Omega,\mathscr{F}_T;C([0,T],\mathbb{R})\right)$ is also a separable Banach space. Besides this $\mathscr{M}^{c,p}(\mathbb{R})$ is a closed subspace of $L^p\left(\Omega,\mathscr{F}_T;C([0,T],\mathbb{R})\right)$ as we proceed previously for $E$, thus $\mathscr{M}^{c,p}(\mathbb{R})$ is separable too.\\Let $\mathscr{M}^{c,p,c}(\mathbb{R})$ (resp.$E_c$) denoted a countable subset dense in $\mathscr{M}^{c,p}(\mathbb{R})$ (resp.$E$). We have
\begin{align}
\left\Vert \sum_{j=1}^m\phi^jx_j-\sum_{j=1}^m\phi^{j,n}x_{j,n}\right\Vert_{\mathscr{M}^{c,p}(E)}&\leq \left\Vert \sum_{j=1}^m(\phi^j-\phi^{j,n})x_{j}\right\Vert_{\mathscr{M}^{c,p}(E)}+\left\Vert \sum_{j=1}^m\phi^{j,n}(x_j-x_{j,n})\right\Vert_{\mathscr{M}^{c,p}(E)}\nonumber\\
&\leq  C_1  \sum_{j=1}^m\left\Vert\phi^j-\phi^{j,n}\right\Vert_{\mathscr{M}^{c,p}(\mathbb{R})}+C_2  \sum_{j=1}^m\left\Vert x_j-x_{j,n}\right\Vert_E,\label{ineq1}
\end{align}
where $x_j\in E$, $\phi^j\in \mathscr{M}^{c,p}(\mathbb{R})$ and $x_{j,n}\in E_c$, $\phi^{j,n}\in\mathscr{M}^{c,p,c}(\mathbb{R})$ for each natural number $n$.
\\For constants that appeared in (\ref{ineq1}), we have set $C_1=\displaystyle\max_{1\leq j\leq m}\left\Vert x_j\right\Vert_E,\quad C_2=\max_{1\leq j\leq m}\sup_{d\geq 1}\left\Vert \phi^{j,d}\right\Vert_{\mathscr{M}^{c,p}(\mathbb{R})}.$ Since any finite Cartesian product of countable sets is also a countable set, and any countably union of countable sets is again a countable set, thus there is an onto mapping from a countable set to $\mathscr{M}^{c,p,c}(\mathbb{R})\otimes E_c$. Consequently the algebraic tensor product $\mathscr{M}^{c,p,c}(\mathbb{R})\otimes E_c$ is a countable dense subset in $\mathscr{M}^{c,p}(\mathbb{R})\otimes E$.\bigskip\\
We will now analyze the class of martingales closed by random variables in $L^p(\Omega,\mathscr{F}_T; E)$.\\
Let $\Phi$ be in $L^p(\Omega,\mathscr{F}_T;E)$. Let $\varepsilon$ be a strictly positive real number.
\\By the fact that $L^p(\Omega,\mathscr{F}_T)\otimes E$ is dense in $L^p(\Omega,\mathscr{F}_T;E)$, we may choose a sequence $\left((\phi_j,x_j)\right)_{1\leq j\leq m}$ in $L^p(\Omega,\mathscr{F}_T)\times E$ satisfying
\begin{eqnarray}\label{ineq2}
\left\Vert \Phi-\sum_{j=1}^m\phi_jx_j\right\Vert_{L^p(\Omega,\mathscr{F}_T;E)}<\varepsilon.
\end{eqnarray}
Since $E$ is separable, then the martingale $\left(\mathbb{E}^{\mathscr{F}_t}(\Phi)\right)_{T\geq t\geq 0}$ possesses a cadlag version $\left(\mathbb{E}_{cad}^{\mathscr{F}_t}(\Phi)\right)_{T\geq t\geq 0}$.
\\Besides this, from the properties of the filtration $\mathbb{F}$ we derive that for each $j$ the real valued martingale $\left(\mathbb{E}^{\mathscr{F}_t}(\phi_j)\right)_{T\geq t\geq 0}$ admits a continuous version denoted by $\left(\mathbb{E}_{con}^{\mathscr{F}_t}(\phi_j)\right)_{T\geq t\geq 0}$.
\\It follows from Doob’s maximal inequality and the contractiveness of the conditional expectation
\begin{eqnarray*}
\mathbb{E}\left(\sup_{T\geq t\geq 0}\left\Vert \mathbb{E}_{cad}^{\mathscr{F}_t}(\Phi)-\sum_{j=1}^m\mathbb{E}_{con}^{\mathscr{F}_t}(\phi_j)x_j\right\Vert_E^p\right)&\lesssim_p & \sup_{T\geq t\geq 0}\mathbb{E}\left(\left\Vert \mathbb{E}_{cad}^{\mathscr{F}_t}(\Phi)-\sum_{j=1}^m\mathbb{E}_{con}^{\mathscr{F}_t}(\phi_j)x_j\right\Vert_E^p\right)\\
&=&\sup_{T\geq t\geq 0}\mathbb{E}\left(\left\Vert \mathbb{E}^{\mathscr{F}_t}(\Phi)-\sum_{j=1}^m\mathbb{E}^{\mathscr{F}_t}(\phi_j)x_j\right\Vert_E^p\right)\\
&\leq &\left\Vert \Phi-\sum_{j=1}^m\phi_jx_j\right\Vert_{L^p(\Omega,\mathscr{F}_T;E)}^p.
\end{eqnarray*}
Combining with $(\ref{ineq2})$ we derive on the one hand that $\left(\mathbb{E}_{cad}^{\mathscr{F}_t}(\Phi)\right)_{T\geq t\geq 0}$ has continuous trajectories, and on the other hand $\mathscr{M}^{c,p}(\mathbb{R})\otimes E$ is dense in the vector space  $\mathscr{M}^{c,p}_{closed}(E)$ consisting of equivalence classes by indistinguishability of continuous versions $\left(\mathbb{E}_{con}^{\mathscr{F}_t}(\Phi)\right)_{T\geq t\geq 0}$, $\Phi\in L^p(\Omega,\mathscr{F}_T;E)$.\\
Finally, from the fact that $\mathscr{M}^{c,p}(E)=\mathscr{M}^{c,p}_{closed}(E)$ the proof is accomplished.
\end{proof}
\begin{thm}\label{thm306}
 Assume that $\left(\Omega,\mathscr{F},\mathbb{P}\right)$ is separable. There exists a sequence $\left((x^n,f^n)\right)_{n\geq 1}$ in $\mathfrak{S}$ such that for every $t>0$ and almost surely $\omega$, we have
$$\left(\int_{[0,t]}^{\mathscr{R}}\mathfrak{S}\, dW_u\right)(\omega)=cl_E\left\lbrace \left(x^n+\int_0^tf^n_udW_u\right)(\omega):n\geq 1\right\rbrace.$$
\end{thm}
\begin{proof}
Let us introduce the set $$\mathbb{I}\left(\mathfrak{S}\right):=\left\lbrace\left(x+\int_0^tf_udW_u\right)_{T\geq t\geq 0},\;(x,f)\in \mathfrak{S}\right\rbrace,$$ constituted by continuous martingales.
According to Lemma \ref{lem301}, $\mathbb{I}\left(\mathfrak{S}\right)$ equipped with the topology induced by $\mathscr{M}^{c,p}(E)$, is also separable. Then, there is a sequence $(f^n)_{n\geq 1}$ of $\mathfrak{S}$ such that 
\begin{eqnarray*}
\mathbb{I}\left(\mathfrak{S}\right)=cl_I\left\lbrace \left(x^n+\int_0^tf_u^ndW_u\right)_{T\geq t\geq 0}:n\geq 1\right\rbrace,
\end{eqnarray*}
where $cl_I$ denotes the closure with respect to the induced topology.\bigskip\\
Let $s$ be in $(0,T]$ and let $h$ be in $\overline{dec}_{\mathscr{F}_s}\left(\pi_{s}\left(\mathbb{I}\left(\mathfrak{S}\right)\right)\right)$.\\
Let $(h^m)_{m\geq 1}$ be a sequence from $dec_{\mathscr{F}_s}\left(\pi_{s}\left(\mathbb{I}\left(\mathfrak{S}\right)\right)\right)$ strongly converging in $L^p\left(\Omega,\mathscr{F}_{s};E\right)$ to $h$.\\
Let for each natural number $m$, $(\Omega_{m,l})_{1\leq l\leq k_m}$ be an $\mathscr{F}_{s}$-measurable partition of $\Omega$, and $(h^{m,l})_{1\leq l\leq k_m}$ be a sequence from $\pi_{s}\left(\mathbb{I}\left(\mathfrak{S}\right)\right)$ with $h^m=\sum_{l=1}^{k_m}\mathds{1}_{\Omega_{m,l}}h^{m,l}$.\bigskip\\
Let for each natural number $m$ and for each $l$, $M^{m,l}$ be from $\mathbb{I}\left(\mathfrak{S}\right)$ such that $\pi_{s}\left(M^{m,l}\right)=h^{m,l}$.\\
Let $\left((x^{\varphi_{m,l}(d)},f^{\varphi_{m,l}(d)})\right)_{d\geq 1}$ be a subsequence of $\left((x^n,f^n)\right)_{n\geq 1}$ such that for each $m$ and $l$
\begin{eqnarray*}
		\left(x^{\varphi_{m,l}(d)}+\int_0^tf_u^{\varphi_{m,l}(d)}dW_u\right)_{T\geq t\geq 0}\xrightarrow[d\to +\infty]{\quad\mathscr{M}^{c,p}(E)\quad}M^{m,l}.
	\end{eqnarray*} 
We have successively
\begin{eqnarray*}
&\mathbb{E}&\left(\left\Vert h-\sum_{l=1}^{k_m}\mathds{1}_{\Omega_{m,l}}\pi_{s}\left(\left(x^{\varphi_{m,l}(d)}+\int_0^tf_u^{\varphi_{m,l}(d)}dW_u\right)_{T\geq t\geq 0}\right)\right\Vert^p\right)\\
&\leq & 2^{p-1}\mathbb{E}\left(\left\Vert h-h^m\right\Vert^p\right)+2^{p-1}\mathbb{E}\left(\sum_{l=1}^{k_m}\mathds{1}_{\Omega_{m,l}}\left\Vert h^{m,l}-\pi_{s}\left(\left(x^{\varphi_{m,l}(d)}+\int_0^tf_u^{\varphi_{m,l}(d)}dW_u\right)_{T\geq t\geq 0}\right)\right\Vert^p\right)\\
&\leq & 2^{p-1}\mathbb{E}\left(\left\Vert h-h^m\right\Vert^p\right)+2^{p-1}\mathbb{E}\left(\max_{1\leq l\leq k_m}\left\Vert h^{m,l}-\pi_{s}\left(\left(x^{\varphi_{m,l}(d)}+\int_0^tf_u^{\varphi_{m,l}(d)}dW_u\right)_{T\geq t\geq 0}\right)\right\Vert^p\right)\\
&\leq & 2^{p-1}\mathbb{E}\left(\left\Vert h-h^m\right\Vert^p\right)+ 2^{p-1}\sum_{l=1}^ {k_m}\mathbb{E}\left(\sup_{T\geq t\geq 0}\left\Vert M_t^{m,l}-x^{\varphi_{m,l}(d)}-\int_0^tf_u^{\varphi_{m,l}(d)}dW_u\right\Vert^p\right).
\end{eqnarray*}
Thus
\begin{eqnarray*}
\overline{dec}_{\mathscr{F}_s}\left(\pi_{s}\left(\mathbb{I}\left(\mathfrak{S}\right)\right)\right)=\overline{dec}_{\mathscr{F}_s}\left(\pi_{s}\left\lbrace x^n+\int_0^tf^n_udW_u:n\geq 1\right\rbrace\right).
\end{eqnarray*}
Subsequently
\begin{eqnarray}\label{assertion1}
\forall s\in (0,T],\;\overline{dec}_{\mathscr{F}_s}\left(\mathbb{I}_{s}\left(\mathfrak{S}\right)\right)=\overline{dec}_{\mathscr{F}_s}\left\lbrace x^n+\int_0^{s}f^n_udW_u:n\geq 1\right\rbrace.
\end{eqnarray}
Define $\varTheta_t(\omega):=cl_E\left\lbrace \left(x^n+\int_0^tf^n_udW_u\right)(\omega):n\geq 1\right\rbrace \textrm{ for all } \omega\in \Omega.$\\
Since the subtrajectory integrals $S^p_{\mathscr{F}_t}(\varTheta_t)$ is a closed decomposable set of $L^p\left(\Omega,\mathscr{F}_t;E\right)$, then
\begin{align*}
S^p_{\mathscr{F}_t}(\varTheta_t)\supset\overline{dec}_{\mathscr{F}_t}\left\lbrace x^n+\int_0^tf^n_udW_u:n\geq 1\right\rbrace.	
\end{align*}
Moreover, Lemma $1.3$ in \cite{Hiai} leads us to the converse inclusion
\begin{align*}
	S^p_{\mathscr{F}_t}(\varTheta_t)\subset\overline{dec}_{\mathscr{F}_t}\left\lbrace x^n+\int_0^tf^n_udW_u:n\geq 1\right\rbrace.
\end{align*}
Therefore, Equality (\ref{assertion1}) and Corollary $1.2$ in \cite{Hiai} together complete the proof.
\end{proof}

\begin{rk}\label{rem3.9}
	By virtue of Lemma $3.3.3$ in \cite{set-stochastic-ch}, the subtrajectory integrals $S^p_{\mathscr{F}_t}\left(\int_{[0,t]}^{\mathscr{R}}\mathfrak{S}\, dW_u\right)$ is a convex subset of the space $L^p(\Omega,\mathscr{F}_t;E)$, and according to Corollary $2.3.3$ in \cite{set-stochastic-ch} we infer that the revised set-valued stochastic integral of $\mathfrak{S}$ is convex-valued since $E$ has the Radon–Nikodym property. On the other hand, Theorem $3.3.5$ in \cite{set-stochastic-ch} asserts that the subtrajectory integrals is a bounded subset of $L^p(\Omega,\mathscr{F}_t;E)$. Thus, Corollary $2.3.4$ in \cite{set-stochastic-ch} yields that  $\int_{[0,t]}^{\mathscr{R}}\mathfrak{S}\, dW_u$  is $p$-integrably bounded, and it takes bounded-valued for almost surely $\omega$.\\
	Moreover, from [\cite{set-stochastic-ch}, Lemma $3.3.2$] and [\cite{Hiai}, Corollary $1.2$] we derive that
	\begin{eqnarray*}
		\int_{[0,t]}^{\mathscr{R}}cl(\mathfrak{S})\, dW_u=\int_{[0,t]}^{\mathscr{R}}\mathfrak{S}\, dW_u\;a.s.\;\omega\in\Omega.
	\end{eqnarray*}
\end{rk}
%\begin{rk}\label{rk302}
%Let us consider the following mapping
%\begin{eqnarray}
%\begin{array}{ccccc}
%\Upsilon & : & E\times \mathcal{L}_{\mathbb{F}}^p([0,T]\times \Omega;E) & \longrightarrow & L^p\left(\Omega,\mathscr{F}_T;E\right)\\
% & & (x,f) & \longmapsto & x+\int_0^Tf_udW_u.\\
%\end{array}\nonumber
%\end{eqnarray}
%Using the martingale representation theorem, one can clearly check that $\Upsilon$ is a linear homeomorphism, which guarantees the reflexivity of the Banach space $E\times\mathcal{L}_{\mathbb{F}}^p([0,T]\times \Omega;E)$. It turns out that a set belongs to $\mathscr{K}_{cwcmpt}\left(E\times\mathcal{L}_{\mathbb{F}}^p([0,T]\times \Omega;E)\right)$ if and only if it is a nonempty, bounded, closed convex subset of the space $E\times\mathcal{L}_{\mathbb{F}}^p([0,T]\times \Omega;E)$.
%\end{rk}
\section{$\gamma$-\;Set-valued Stochastic Integral}
Let $E$ be a separable UMD space, and assume that $\mathbb{F}$ is the natural filtration generated by a standard Brownian motion $W=\left(W_t\right)_{t\in [0,T]}$ augmented by all $\mathbb{P}$-null sets of $\mathscr{F}$.\\
Let $H:[0,T]\times \Omega\longrightarrow \mathscr{K}(E)$ be a measurable and adapted set-valued function.\\
Let us consider the subset $\Lambda_{p,H}^{0,T}$ defined as follows:
\begin{eqnarray}\label{formula2}
\Lambda_{p,H}^{0,T}:=\Big\{ f\in\mathcal{L}_{\mathbb{F}}^p([0,T]\times \Omega;E):f(v,\omega)\in H(v,\omega)\;a.e.(v,\omega)\in [0,T]\times \Omega\Big\},
\end{eqnarray}
where $\mathcal{L}_{\mathbb{F}}^p([0,T]\times \Omega;E)$ is the Banach space stated in Proposition \ref{prop301}.\bigskip\\
%To construct the set-valued integral discussed in this subsection, we need to apply the functional $\mathbb{J}_{s,t}$ on each $f\in \Lambda_{p,H}^{0,T}$ in order to obtain the function $\omega\longmapsto\int_s^tf(u,\omega)du$, which, according to Proposition \ref{prop247}, is an element of $L^p(\Omega,\mathscr{F}_t;E)$. For almost surely $\omega$, the notation $\int_s^t f(u,\omega)du$ is used for the sake of readability, as described in Proposition \ref{prop245}, and is considered to solve the BSDE (\ref{bdeq1}) within the framework of UMD spaces.
Based on Proposition \ref{prop247} and Theorem $3.8$ in [\cite{papa}, Chapter $2$], we construct the following notion.
\begin{defi}\label{defi302}
Assume that $\Lambda_{p,H}^{0,T}$ is nonempty. Let $s,t$ be two real numbers in $[0,T]$ such that $s<t$, and consider the following subset of $L^p(\Omega,\mathscr{F}_t;E)$:
$$\mathbb{J}_{s,t}(\Lambda_{p,H}^{0,T}):=\left\lbrace \int_s^tf(u,.)du:\;f\in \Lambda_{p,H}^{0,T}\right\rbrace.$$ 
The $\gamma$-set-valued stochastic integral of $H$ over $[s,t]$ is the $\mathscr{F}_t$-measurable set-valued random variable, denoted by $\int_s^t H_u du$, which satisfies
\begin{eqnarray}\label{assertion8}
S^p_{\mathscr{F}_t}\left(\int_s^t H_u du\right)=\overline{dec}_{\mathscr{F}_t}\left(\mathbb{J}_{s,t}(\Lambda_{p,H}^{0,T})\right).	
\end{eqnarray} 
\\It is convenient to put $\displaystyle\int_0^0H_u du=\left\lbrace 0_E\right\rbrace=\int_T^TH_u du$ for $a.s.\;\omega \in\Omega$.
\end{defi}
\begin{com}
	Let $p=2$.
	By [\cite{jinping}, Proposition $2.5$], the set-valued $H$ is $\Sigma_{\mathbb{F}}$-measurable. One may then inquire about the relationship between the subtrajectory integrals $S^2_{\Sigma_{\mathbb{F}}}(H)$ and the set $\Lambda_{2,H}^{0,T}$, that ensures compatibility in the context of $L^2$-stochastically integrable processes.\\ For a Hilbert space $E$, we have the isometry $\gamma(0,T;E)= L^2([0,T];E)$ according to [\cite{neerven2}, Proposition 9.2.9]. Thus, $\Lambda_{2,H}^{0,T}=S^2_{\Sigma_{\mathbb{F}}}(H)$. Moreover, [\cite{neerven2}, Theorem 7.3.1] and Proposition \ref{prop246} assert that the random variable constructed in \eqref{assertion8} coincides with the Lebesgue set-valued stochastic integral defined by \eqref{lebint} for the spaces $E=\R^m$ and the $2$-integrably boundedness condition on $H$. 
%	The reader may also refer to  [\cite{kis2012}, Section 3] for the set-valued stochastic integral of an m-dimensional Itô integrable set-valued process with respect to the Lebesgue measure $dt$. 
%	Let $p=2$. Due to Proposition $2.5$ in \cite{jinping}, $H$ is $\Sigma_{\mathbb{F}}$-measurable, and thanks to [\cite{intermeas}, Theorem $6.6.8$], $H$ admits a Casting representation $(f^n)_{n\geq 1}$ of $\Sigma_{\mathbb{F}}$-measurable selections.
%	\\By supplying $E$ with the property type $2$ (resp. cotype $2$) we get $\Lambda_{2,H}^{0,T}\supset S^2_{\Sigma_{\mathbb{F}}}(H)$ \big(resp. $\Lambda_{2,H}^{0,T}\subset S^2_{\Sigma_{\mathbb{F}}}(H)$\big), and if $E$ is a Hilbert space we have $\Lambda_{2,H}^{0,T}=S^2_{\Sigma_{\mathbb{F}}}(H)$ due to Proposition $9.2.9$ in \cite{neerven2}.
\end{com}
\noindent Next, we show the localization property of the $\gamma$-set-valued stochastic integral.
\begin{lem}\label{lem302}
Assume that $\Lambda_{p,H}^{0,T}$ is nonempty. For every $s,t\in [0,T]$ such that $s<t$, the following statement holds
$$\int_s^t H_u du=\int_0^T \mathds{1}_{[s,t]} H_u du.$$
\end{lem}
\begin{proof}
Taking into account Remark \ref{rk241} and Proposition \ref{prop245}, we conclude that $\Lambda_{p,\mathds{1}_{[s,t]} H}^{0,T}$ is nonempty, and we have $$\mathbb{J}_{s,t}(\Lambda_{p,H}^{0,T})=\mathbb{J}_{0,T}(\Lambda_{p,\mathds{1}_{[s,t]} H}^{0,T}).$$
Hence $$\mathbb{J}_{s,t}(\Lambda_{p,H}^{0,T})\subset\overline{dec}_{\mathscr{F}_t}\left(\mathbb{J}_{0,T}(\Lambda_{p,\mathds{1}_{[s,t]} H}^{0,T})\right).$$
Consequently
\begin{eqnarray*}
\overline{dec}_{\mathscr{F}_t}\left(\mathbb{J}_{s,t}(\Lambda_{p,H}^{0,T})\right)\subset\overline{dec}_{\mathscr{F}_T}\left(\mathbb{J}_{0,T}(\Lambda_{p,\mathds{1}_{[s,t]} H}^{0,T})\right).	
\end{eqnarray*} 
On the other hand, we have
\begin{eqnarray*}
S^p_{\mathscr{F}_T}\left(\int_0^T \mathds{1}_{[s,t]}H_udu\right)=\overline{dec}_{\mathscr{F}_T}(\mathbb{J}_{s,t}(\Lambda_{p,H}^{0,T}))\subset \overline{dec}_{\mathscr{F}_T}\left(S^p_{\mathscr{F}_t}\left(\int_s^tH_udu\right)\right)\subset S^p_{\mathscr{F}_T}\left(\int_s^tH_udu\right).	
\end{eqnarray*}
Therefore, 
$$S^p_{\mathscr{F}_T}\left(\int_0^T \mathds{1}_{[s,t]}H_udu\right)=S^p_{\mathscr{F}_T}\left(\int_s^tH_udu\right).$$
Finally, it suffices to apply Corollary $1.2$ in \cite{Hiai}.
\end{proof}
\noindent To explore additional results, we need the following assumption :  \begin{align}\label{assertion6}  \exists \nu \in L^p(\Omega, \mathscr{F}_T),\,\,\forall f \in \Lambda_{p,H}^{0,T},\quad\left\Vert f(\cdot,\omega) \right\Vert_{\gamma(0,T;E)} \leq \nu, \quad \text{a.s.}\,\, \omega \in \Omega\tag{$\mathscr{H}_0$}. \end{align}
We now establish, by the following theorem, that the $\gamma$-set-valued stochastic integral satisfies the Chasles relation, thereby contributing to the analysis used to solve the SVBSDE (\ref{sveq1}).
\begin{thm}\label{thm311}
Assume that $\Lambda_{p,H}^{0,T}$ is nonempty and satisfies assumption (\ref{assertion6}). For every $t\in (0,T)$, we have
\begin{eqnarray*}
\int_0^TH_u du &=& \int_0^t H_u du+\int_t^TH_u du\quad a.s.\;\omega\in\Omega;\\
\int_t^TH_u du &=& \int_0^TH_u du \circleddash \int_0^tH_u du\quad a.s.\;\omega\in\Omega.
\end{eqnarray*}
\end{thm}
\begin{proof}
Let $t\in (0,T)$. By virtue of Lemma $2.7$ in \cite{neerven3} we have for every $f\in \Lambda_{p,H}^{0,T}$:
\begin{eqnarray*}
\left\Vert \int_{s_1}^{s_2} f(u,\omega)du\right\Vert_E &\leq & (s_2-s_1)^{\frac{1}{2}} \nu\quad a.s.\omega\in\Omega,
\end{eqnarray*}
where $s_1,s_2\in [0,T],\;s_1<s_2.$ Thus, $\mathbb{J}_{0,t}(\Lambda_{p,H}^{0,T}),\;\mathbb{J}_{t,T}(\Lambda_{p,H}^{0,T})$ and $\mathbb{J}_{0,T}(\Lambda_{p,H}^{0,T})$ are $p$-integrably bounded. 
Besides this, $\overline{dec}_{\mathscr{F}_T}\left(\mathbb{J}_{0,T}(\Lambda_{p,\mathds{1}_{[0,t]}H}^{0,T})\right)$ and $\overline{dec}_{\mathscr{F}_T}\left(\mathbb{J}_{0,T}(\Lambda_{p,\mathds{1}_{[t,T]}H}^{0,T})\right)$ are weakly closed since the closed decomposable hull of a nonempty convex set contained in $L^p\left(\Omega,\mathscr{F}_T;E\right)$ is again a convex set. Further, $L^p\left(\Omega,\mathscr{F}_T;E\right)$ is reflexive and the decomposable hull of a $p$-integrably bounded subset of $L^p\left(\Omega,\mathscr{F}_T;E\right)$ is in fact a bounded subset of this Bochner space.\\Then, $\overline{dec}_{\mathscr{F}_T}\left(\mathbb{J}_{0,T}(\Lambda_{p,\mathds{1}_{[t,T]}H}^{0,T})\right)$ is weakly compact.

Next, thanks to Proposition \ref{prop245} and [\cite{Hiai}, Theorem $1.4$] we have
\begin{eqnarray}
S^p_{\mathscr{F}_{T}}\Big(\int_{0}^{T}H_u du\Big)&=&\overline{dec}_{\mathscr{F}_{T}}\left(\mathbb{J}_{0,t}(\Lambda_{p,H}^{0,T})+\mathbb{J}_{t,T}(\Lambda_{p,H}^{0,T})\right)\nonumber\\
&=&\overline{dec}_{\mathscr{F}_{T}}\left(\mathbb{J}_{0,T}(\Lambda_{p,\mathds{1}_{[0,t]}H}^{0,T})+\mathbb{J}_{0,T}(\Lambda_{p,\mathds{1}_{[t,T]}H}^{0,T})\right)\nonumber\\
&=& S^p_{\mathscr{F}_{T}}\Big(\int_{0}^{T}\mathds{1}_{[0,t]}H_u du\Big)+S^p_{\mathscr{F}_{T}}\Big(\int_{0}^{T}\mathds{1}_{[t,T]}H_udu\Big)\nonumber\\
&=& cl_{L^p(\Omega,\mathscr{F}_T;E)}\Big(S^p_{\mathscr{F}_{T}}\Big(\int_{0}^{T}\mathds{1}_{[0,t]}H_u du\Big)+S^p_{\mathscr{F}_{T}}\Big(\int_{0}^{T}\mathds{1}_{[t,T]}H_u du\Big)\Big)\nonumber\\
&=& S^p_{\mathscr{F}_{T}}\Big(cl_E\Big(\int_{0}^{T}\mathds{1}_{[0,t]}H_u du+\int_{0}^{T}\mathds{1}_{[t,T]}H_udu\Big)\Big).\label{assertion5}
\end{eqnarray}
Now, by arguing similarly as in Remark \ref{rem3.9}, we obtain that $\int_0^T H_u du$, $\int_0^t H_u du$ and $\int_t^T H_u du$ have convex weakly compact-valued for almost surely $\omega$. Therefore, by coming back to (\ref{assertion5}) we get the desired result.
\end{proof}

Here, we provide a simple example of set-valued functions that satisfy some important properties.
\begin{examp}
Let $\tilde{H}:[0,T]\times \Omega\longrightarrow \mathscr{K}(\mathbb{R})$ be a $\mathscr{B}([0,T])\otimes\mathscr{F}$-measurable $\mathbb{F}$-adapted set-valued function such that $\mathbb{E}\left(\left\Vert \vert\tilde{H}\vert_{\varrho'}\right\Vert_{L^2([0,T])}^p\right)$ is finite, where $\varrho'$ is the Hausdorff distance on $\mathscr{K}(\mathbb{R})$.\\
Let $x$ be a fixed element of $E$, and let's put $H(t,\omega)=\tilde{H}(t,\omega)x$ for every $(t,\omega)$ in $[0,T]\times \Omega$.\\
Let $f$ be a $\Sigma_{\mathbb{F}}$-measurable selection of $\tilde{H}$, and $A$ be an $\mathscr{F}_T$-measurable set.\\
Thanks to Corollary $1.1.2$ in \cite{neerven1}, the function $(t,\omega)\mapsto f(t,\omega)x$ is $\Sigma_{\mathbb{F}}$-measurable, and by virtue of Example \ref{examp1} we have:
$$\mathbb{E}\left(\mathds{1}_A\left\Vert f\otimes x\right\Vert_{\gamma(0,T;E)}^p\right)=\mathbb{E}\left(\mathds{1}_A\left\Vert f\right\Vert_{L^2([0,T])}^p\left\Vert x\right\Vert_{E}^p\right)\leq\left\Vert x\right\Vert_{E}^p\mathbb{E}\left(\mathds{1}_A\left\Vert \vert\tilde{H}\vert_{\varrho'}\right\Vert_{L^2([0,T])}^p\right).$$
Thus, $\Lambda_{p,H}^{0,T}$ is nonempty and the assumption (\ref{assertion6}) is fulfilled.
\end{examp}

We show by the following lemma that we have a form of commutation when applying $\mathbb{J}_{s,t}$ and the closure.
\begin{lem}\label{lem400}
Assume that $\Lambda_{p,H}^{0,T}$ is nonempty and satisfies assumption (\ref{assertion6}), and that $H$ is convex valued. Let $s,t$ be in $[0,T]$ with $s<t$. The following assertion holds:
\begin{eqnarray}\label{eq401}
cl_{L^p}\left(\mathbb{J}_{s,t}(\Lambda_{p,H}^{0,T})\right)=\mathbb{J}_{s,t}\left(cl_{\mathcal{L}_{\mathbb{F}}^p}\left(\Lambda_{p,H}^{0,T}\right)\right),	
\end{eqnarray}
where "$cl_{\mathcal{L}_{\mathbb{F}}^p}$" (resp. "$cl_{L^p}$") denotes the closure in $\mathcal{L}_{\mathbb{F}}^p([0,T]\times \Omega;E)$ (resp. $L^p\left(\Omega,\mathscr{F}_t;E\right)$). 	
\end{lem}
\begin{proof}
Let $(f^n))_{n \geq 1}$ be a sequence in $\Lambda_{p,H}^{0,T}$ such that $\left(\mathbb{J}_{s,t}( f^n)\right)_{n \geq 1}$ converges in $L^p(\Omega, \mathscr{F}_t; E)$.
The linear mapping:
\begin{eqnarray}
	\begin{array}{ccccc}
		&  &  \mathcal{L}_{\mathbb{F}}^p([0,T]\times \Omega;E) & \longrightarrow & L_0^p\left(\Omega,\mathscr{F}_T;E\right)\\
		& & g & \longmapsto & \displaystyle\int_0^Tg_udW_u,\\
	\end{array}\label{map2}
\end{eqnarray}
establishes an isomorphism from $ \mathcal{L}^p_{\mathbb{F}}\left([0, T] \times \Omega; E\right)$ onto $L_0^p(\Omega, \mathscr{F}_T; E)$, and $L_0^p(\Omega, \mathscr{F}_T; E)$ is a reflexive Banach space according to Theorem $1.11.16$ in \cite{reflexive} . Thus, $ \mathcal{L}^p_{\mathbb{F}}\left([0, T] \times \Omega; E\right)$ is itself reflexive. Subsequently, there is a weakly convergent subsequence $(f^{\varphi(n)})_{n \geq 1}$ in $\mathcal{L}^p_{\mathbb{F}}\left([0, T] \times \Omega; E\right)$ to an element $f$, and we have  $f \in cl_{\mathcal{L}_{\mathbb{F}}^p}(\Lambda_{p,H}^{0,T})$ since $\Lambda_{p,H}^{0,T}$ is convex.\\
By virtue of Proposition \ref{prop245}, we derive that
 $\left(\mathbb{J}_{s,t}( f^{\varphi(n)})\right)_{n \geq 1}$ is weakly convergent in $L^p(\Omega, \mathscr{F}_t; E)$ to $\mathbb{J}_{s,t}( f)$.
Thus, $cl_{L^p}\left(\mathbb{J}_{s,t}(\Lambda_{p,H}^{0,T})\right)$ is contained in $\mathbb{J}_{s,t}\left(cl_{\mathcal{L}_{\mathbb{F}}^p}\left(\Lambda_{p,H}^{0,T}\right)\right)$.\\We obtain the reverse inclusion by continuity, and therefore, equality (\ref{eq401}) is achieved.
\end{proof}
Now, we can state the following representation theorem.
\begin{thm}
	Assume that $\Lambda_{p,H}^{0,T}$ is a nonempty closed subset of $\mathcal{L}^p_{\mathbb{F}}\left([0, T] \times \Omega; E\right)$, satisfying assumption (\ref{assertion6}), $H$ is convex valued, and  $(\Omega,\mathscr{F},\mathbb{P})$ is separable. Then, there exists a sequence $(f^n)_{n\geq 1}$ contained in $\Lambda_{p,H}^{0,T}$ such that for every $s,t\in [0,T],\;s<t$, we have
	\begin{eqnarray}\label{represent}
	\left(\int_s^t H_udu\right)(\omega)=cl_E\left\lbrace \int_s^tf^n(u,\omega)du:n\geq 1\right\rbrace\;a.s.\,\omega\in\Omega.	
	\end{eqnarray}
\end{thm}
\begin{proof}
Since $\mathcal{L}^p_{\mathbb{F}}\left([0, T] \times \Omega; E\right)$ is separable,  then we can assign a sequence $(f^n)_{n\geq 1}$ in this space such that $\Lambda_{p,H}^{0,T}=cl_{\mathcal{L}_{\mathbb{F}}^p}\left\lbrace f^n:n\geq 1\right\rbrace$. It follows that $\mathbb{J}_{s,t}(\Lambda_{p,H}^{0,T})=cl_{L^p(\Omega, \mathscr{F}_t; E)}\left\lbrace\int_s^tf^n(u,.)du:n\geq 1\right\rbrace$ according to Lemma \ref{lem400}. Next, Lemma $3.3.2$ in \cite{set-stochastic-ch} asserts that
$$S^p_{\mathscr{F}_t}\left(\int_s^t H_u du\right)=\overline{dec}_{\mathscr{F}_t}\left(\left\lbrace\int_s^tf^n(u,.)du:n\geq 1\right\rbrace\right).$$
By considering the following $\mathscr{F}_t$-measurable set-valued function
$$\varTheta_{s,t}(\omega):=cl_E\left\lbrace\int_s^tf^n(u,\omega)du:n\geq 1\right\rbrace\quad\textrm{a.s. }\omega,$$ 
we get the aim by virtue of Lemma $1.3$ in \cite{Hiai}.
\end{proof}
\begin{rk}
The representation (\ref{represent}) remains valid if we modify the assumptions by removing the condition that '$\Lambda_{p,H}^{0,T}$ is closed' and by endowing the space $E$ with the geometric property cotype $2$. Indeed, if $(g^n)_{n\geq 1}$ is a sequence from within $\Lambda_{p,H}^{0,T}$ conveging to g in $\mathcal{L}^p_{\mathbb{F}}\left([0, T] \times \Omega; E\right)$, we derive that it converges in $L^p(\Omega;L^2([0,T];E))$ according to Theorem $9.2.11$ in \cite{neerven2}. Subsequently, it converges in measure on $[0,T]\times\Omega$ to $g$. Taking into account that $H$ is closed valued, we conclude that $g\in \Lambda_{p,H}^{0,T}$.
\end{rk}

To institute a strengthened foundation for the 
$\gamma$-set-valued stochastic integral, we introduce another concept by generalizing the Lebesgue set-valued stochastic integral defined by (\ref{lebint}).
\begin{defi}
	Let $p=2$ and $H$ be $2$-integrably bounded. Let $s,t\in [0,T]$ such that $s<t$.
	\\The Bochner set-valued integral of $H$ over $[s,t]$ is the set-valued random variable defined by 
	\begin{eqnarray}\label{bocint}
		S_{\mathscr{F}_t}^2\left((\mathsf{b}_{\mathbb{F}})\int_s^t H_udu\right)=\overline{dec}_{\mathscr{F}_t}\left(\left\lbrace \int_{[s,t]}^{Bochner}g(u,.)du:g\in S_{\Sigma_{\mathbb{F}}}^2(H)\right\rbrace\right).
	\end{eqnarray}
\end{defi}
\noindent Next, the following result aims to establish a relationship between the $\gamma$-set-valued stochastic integral and the Bochner set-valued integral under suitable conditions. It is based on comparing the subtrajectory integrals and the set defined by (\ref{formula2}).
\begin{thm}\label{thm403}
	Let $p=2$ and $H$ be $2$-integrably bounded. Let $s,t$ be in $[0,T]$ such that $s<t$.
	\\The following assertions hold:
	\begin{itemize}
		\item[$(1)$] If the space $E$ has type $2$, we have: 
		\begin{eqnarray*}
			(\mathsf{b}_{\mathbb{F}})\int_s^t H_udu\subset \int_s^t H_udu\quad \textrm{a.s.}\; \omega\in\Omega.
		\end{eqnarray*}
		\item[$(2)$] If $\Lambda_{2,H}^{0,T}$ is nonempty and the space $E$ has cotype $2$, we have:
		\begin{eqnarray*}
			(\mathsf{b}_{\mathbb{F}})\int_s^t H_udu\supset \int_s^t H_udu\quad \textrm{a.s.}\; \omega\in\Omega.
		\end{eqnarray*}
	\end{itemize}	
\end{thm}
\begin{proof}
\begin{itemize}
	\item [(1)] Assume that $E$ has type $2$, and $S_{\Sigma_{\mathbb{F}}}^2(H)$ is nonempty.\\
	Let $\varphi:[0,T]\times\Omega\longrightarrow E$ be a $\Sigma_{\mathbb{F}}$-measurable function in $S_{\Sigma_{\mathbb{F}}}^2(H)$. According to Fubini theorem and the separability of the space $L^2([0,T])$, there exists a negligible set $\mathscr{N}\in \mathscr{F}_0$ such that for any $\omega\in \Omega\setminus\mathscr{N}$, the mapping $$X(\omega):h\in L^2([0,T])\longmapsto \int^{Bochner}_{[0,T]}h(u)\varphi(u,\omega)du$$
	is well defined, and the function $\varphi(.,\omega)$ belongs to $L^2([0,T];E)$. It follows by virtue of [\cite{neerven2}, Theorem $9.2.10$] that $X(\omega)$ is an element of $\gamma(0,T;E)$ for each $\omega\in \Omega\setminus\mathscr{N}$.\\
	Fix $\mathbb{T}$ in $\gamma(0,T;E)$ and define $\tilde{X}:=\mathds{1}_{\mathscr{N}}\mathbb{T}+\mathds{1}_{\Omega\setminus\mathscr{N}}X$.
	\\By applying once more Fubini theorem, we derive that $\omega\longmapsto \tilde{X}(\omega)h$ is a $\mathbb{P}$-version of an $\mathscr{F}_T$-measurable function for each $h\in L^2([0,T])$. Thus, Lemma $2.7$ in \cite{neerven3} asserts that $\tilde{X}$ is a random variable. On the other hand, we have for any $\omega\in \Omega\setminus\mathscr{N}$
	$$\left\Vert \tilde{X}(\omega) \right\Vert_{\gamma(0,T;E)}\leq \varsigma_{2,E}^{\gamma}\left\Vert \varphi(.,\omega)\right\Vert_{L^2([0,T];E)},$$
	where $\varsigma_{2,E}^{\gamma}$ is the Gaussian type $2$ constant (\cite{neerven2}, Definition $7.1.17$). Then $\tilde{X}\in L^2(\Omega;\gamma(0,T;E))$.\\ Besides this, $\varphi$ defines $\tilde{X}$ according to Lemma $2.7$ in \cite{neerven3}. Therefore, $\varphi$ is $L^2$-stochastically integrable, and $S_{\Sigma_{\mathbb{F}}}^2(H)$ is contained in $\Lambda_{2,H}^{0,T}$. Next, Proposition \ref{prop246} guarantees that $$\left\lbrace \int_{[s,t]}^{Bochner}g(u,.)du:g\in S_{\Sigma_{\mathbb{F}}}^2(H)\right\rbrace\subset \mathbb{J}_{s,t}(\Lambda_{p,H}^{0,T}).$$
	Thus, we achieve the first aim by applying Lemma $2.3$ in \cite{jinping0}.
	\item [$(2)$]Assume that $E$ has cotype $2$, and $\Lambda_{2,H}^{0,T}$ is nonempty.\\
	Let $\varphi:[0,T]\times\Omega\longrightarrow E$ be a $\Sigma_{\mathbb{F}}$-measurable function in $\Lambda_{2,H}^{0,T}$.\\
	Since $\varphi$ defines a random variable $X$ in $L^2(\Omega;\gamma(0,T;E))$, there exists a negligible set $\mathscr{N}$ such that $\varphi(.,\omega)$ defines $X(\omega)$ for any $\omega\in\Omega\setminus \mathscr{N}$.\\
	Moreover, by virtue of Theorem $9.2.11$ in \cite{neerven2} we have $$\forall \omega\in \Omega\setminus \mathscr{N},\quad\left\Vert \varphi(.,\omega)\right\Vert_{L^2([0,T];E)}\leq\upsilon_{2,E}^{\gamma}\left\Vert X(\omega) \right\Vert_{\gamma(0,T;E)},$$
	where $\upsilon_{2,E}^{\gamma}$ is the Gaussian cotype $2$ constant (\cite{neerven2}, Definition $7.1.17$). Therefore, $\varphi$ belongs to $S_{\Sigma_{\mathbb{F}}}^2(H)$, and $\Lambda_{2,H}^{0,T}$ is contained in $S_{\Sigma_{\mathbb{F}}}^2(H)$. It follows by virtue of Proposition \ref{prop246} that $$\mathbb{J}_{s,t}(\Lambda_{p,H}^{0,T})\subset\left\lbrace \int_{[s,t]}^{Bochner}g(u,.)du:g\in S_{\Sigma_{\mathbb{F}}}^2(H)\right\rbrace.$$
\end{itemize}
Finally, Lemma $2.3$ in \cite{jinping0} gives the desired result.	
\end{proof}
Now, by Kwapień’s Theorem 7.3.1 in \cite{neerven2}, we emphasize the particularity of Hilbert spaces and provide a straightforward result from Theorem \ref{thm403}.
\begin{cor}
Let $p=2$ and $E$ be a Hilbert space. Let $s,t\in [0,T]$ such that $s<t$.\\
The set-valued $H$ satisfies assumption (\ref{assertion6}) if and only if $H$ is $2$-integrably bounded. In this case, we have
	\begin{eqnarray*}
		(\mathsf{b}_{\mathbb{F}})\int_s^t H_udu= \int_s^t H_udu\quad \textrm{a.s.}\; \omega\in\Omega.
	\end{eqnarray*}	
\end{cor}	
%\begin{rk}\label{rk301}
%	 Due to Theorem $2.2$ in \cite{Hiai} we have
%	\begin{eqnarray*}
%		\int_{\Omega}\nu^2d\mathbb{P}\geq \sup_{f\in \Lambda_{2,H}^{0,T}}\int_{[0,T]\times \Omega}\left\Vert f\right\Vert^2_E dud\mathbb{P}=\int_{[0,T]\times \Omega}\sup_{n\geq 1}\left\Vert f^n\right\Vert^2_E dud\mathbb{P}=\int_{[0,T]\times \Omega}\left\vert H\right\vert^2_{\varrho}dud\mathbb{P}.
%	\end{eqnarray*}
%\end{rk}
\section{Martingale and Casting Representation Theorems for Set-valued Martingales in UMD Banach Spaces}
We adopt certain constraints to represent a set-valued martingale as a revised set-valued stochastic integral, and establish a Casting representation of this set-valued martingale through a sequence of martingale selections.  After that, under suitable conditions, we show the possibility of representing the Hukuhara difference of two set-valued martingales as a revised set-valued stochastic integral.\\

Let $E$ be a separable UMD space and assume that $\mathbb{F}$ is the natural filtration generated by a standard Brownian motion $W=\left(W_t\right)_{t\in [0,T]}$ augmented by all $\mathbb{P}$-null sets of $\mathscr{F}$.\bigskip\\
Let $\Psi:[0,T] \times \Omega \to \mathscr{K}_c(E)$ be an $\mathbb{F}$-adapted set-valued process such that $\Psi_T$ is $p$-integrably bounded. 

\subsection{Martingale Representation}\label{subsec5.1}
\begin{lem}\label{lem401}
If $\Psi$ is a set-valued martingale, then the set $\mathscr{M}^p_{sel}(\Psi)$ of all $\mathbb{F}$-martingale selections of $\Psi$ is nonempty and forms a closed subset of $\mathscr{M}^{c,p}(E)$.
\end{lem}
\begin{proof}
Let $m$ be the least natural number greater than or equal to $T$. We introduce $\Psi^d=\left(\Psi_{k\wedge T}\right)_{m\geq k\geq 0}$ which is a discretization of $\Psi$. By virtue of Theorem $3.2$ in \cite{hess} we infer that $\mathscr{M}^p_{sel}\left(\Psi^d\right)$ is nonempty and for an element $\left(\psi^d_{k\wedge T}\right)_{m\geq k\geq 0}$ of $\mathscr{M}^p_{sel}\left(\Psi^d\right)$, the process $(\psi_t)_{T\geq t\geq 0}$ constructed in the following manner
\begin{eqnarray}
\psi_t=\sum_{k=1}^m\mathds{1}_{[k-1,k\wedge T)}(t)\mathbb{E}_{con}^{\mathscr{F}_t}\left(\psi^d_{k\wedge T}\right)+\mathds{1}_{\left\lbrace T\right\rbrace}(t)\psi^d_{T}\quad\textrm{in }L^p\left(\Omega;E\right)\;\textrm{ for each }t,
\end{eqnarray}
ensures us the non-emptiness of the set $\mathscr{M}^p_{sel}\left(\Psi\right)$. Indeed, for any $t$, there is $k$ such that
$$\psi_t=\mathbb{E}_{con}^{\mathscr{F}_t}\left(\psi^d_{k\wedge T}\right)\in S^p_{\mathscr{F}_t}\left(\mathbb{E}^{\mathscr{F}_t}\left(\Psi_{k\wedge T}^d \right)\right)=S^p_{\mathscr{F}_t}\left(\Psi_t\right).$$
On the other hand, it is obvious that any martingale selection of $\Psi$ admits a continuous version, which will be the associated representative. Thus, we shall consider $\mathscr{M}^p_{sel}\left(\Psi\right)$ as a subset of $\mathscr{M}^{c,p}\left(E\right)$.\\
Let $\left(\psi^k\right)_{k\geq 1}$ be a sequence in $\mathscr{M}^p_{sel}\left(\Psi\right)$ that converges to $\psi$ in $\mathscr{M}^{c,p}\left(E\right)$. We have
\begin{eqnarray*}
d_{L^p\left(\Omega,\mathscr{F}_t;E\right)}^p\left(\psi_t,S^p_{\mathscr{F}_t}\left(\Psi_t\right)\right)&\leq & 2^{p-1}\left(\mathbb{E}\left(\left\Vert \psi_t-\psi_t^k\right\Vert^p\right)+d_{L^p\left(\Omega,\mathscr{F}_t;E\right)}^p\left(\psi^k_t,S^p_{\mathscr{F}_t}\left(\Psi_t\right)\right)\right)\\
&\leq & 2^{p-1} \left\Vert \psi-\psi^k\right\Vert_{\mathscr{M}^{c,p}\left(E\right)}^p,
\end{eqnarray*}
for any $t\in [0,T]$, where $d_{L^p\left(\Omega,\mathscr{F}_t;E\right)}$ denotes the distance from an element to a non-empty subset of $L^p\left(\Omega,\mathscr{F}_t;E\right)$. Thus, the aim is accomplished.
\end{proof}
\begin{lem}\label{lem402}
If $\Psi$ is a set-valued martingale, then for every $t\in [0,T]$ we have
\begin{align*} \overline{dec}_{\mathscr{F}_t}\left(\pi_t\left(\mathscr{M}^p_{sel}\left(\Psi\right)\right)\right)=S^p_{\mathscr{F}_t}\left(\Psi_t\right).
	\end{align*}
\end{lem}
\begin{proof}
Let $t\in [0,T]$. We have
\begin{eqnarray}\label{assertion7}
\pi_t\left(\mathscr{M}^p_{sel}\left(\Psi\right)\right)\subset \left\lbrace \mathbb{E}_{con}^{\mathscr{F}_t}\left(h\right):h\in S^p_{\mathscr{F}_{T}}\left(\Psi_{T}\right)\right\rbrace
\end{eqnarray}
in $L^p(\Omega,\mathscr{F}_t;E)$.\\Since the set on the right-hand side of (\ref{assertion7}) is a decomposable subset of $L^p\left(\Omega,\mathscr{F}_t;E\right)$, then
\begin{eqnarray*}
\overline{dec}_{\mathscr{F}_t}\left(\pi_t\left(\mathscr{M}^p_{sel}\left(\Psi\right)\right)\right)&\subset & cl_{L^p\left(\Omega,\mathscr{F}_t;E\right)}\Big(\left\lbrace \mathbb{E}^{\mathscr{F}_t}_{con}\left(h\right):h\in S^p_{\mathscr{F}_T}\left(\Psi_T\right)\right\rbrace\Big)\\&= & S^p_{\mathscr{F}_t}\left(\mathbb{E}^{\mathscr{F}_t}\left(\Psi_T\right)\right)= S^p_{\mathscr{F}_t}\left(\Psi_t\right)\\
&\subset &cl_{L^p\left(\Omega,\mathscr{F}_t;E\right)}\Big(\left\lbrace \pi_t(\psi):\psi\in \mathscr{M}^p_{sel}\left(\Psi\right)\right\rbrace\Big)\\
& \subset & cl_{L^p\left(\Omega,\mathscr{F}_t;E\right)}\Big(dec_{\mathscr{F}_t}\left\lbrace \pi_t(\psi):\psi\in \mathscr{M}^p_{sel}\left(\Psi\right)\right\rbrace\Big)\\
&=& \overline{dec}_{\mathscr{F}_t}\left(\pi_t\left(\mathscr{M}^p_{sel}\left(\Psi\right)\right)\right).
\end{eqnarray*}
Finally, the desired formula is proved.
\end{proof}
\begin{lem}\label{lem403}
If $\Psi$ is a set-valued martingale, then the following set
\begin{eqnarray}\label{ineq10}
\mathfrak{S}_{\Psi}:=\left\lbrace (x^{\psi},f^{\psi})\in E\times\mathcal{L}_{\mathbb{F}}^p([0,T]\times \Omega;E):\psi\in \mathscr{M}_{sel}^p(\Psi),\;\psi_.=x^{\psi}+\int_0^{.}f_u^{\psi}dW_u\right\rbrace
\end{eqnarray}
is a nonempty closed convex bounded subset of $E\times\mathcal{L}_{\mathbb{F}}^p([0,T]\times \Omega;E)$.
\end{lem}
\begin{proof}
Let $\upsilon\in L^p(\Omega,\mathscr{F}_T)$ such that $\left\vert \Psi_T\right\vert_{\varrho}\leq \upsilon$ for almost surely $\omega$.\\
First, the non-emptiness of the subset $\mathfrak{S}_{\Psi}$ follows from the non-emptiness of $\mathscr{M}_{sel}^p\left(\Psi\right)$ and the Brownian martingale representation theorem [\cite{BSEE.B}, $(3.2)$].\\
Let $\left((x^n,g^n)\right)_{n\geq 1}$ be a sequence in $\mathfrak{S}_{\Psi}$ converging to $(x,g)\in E\times\mathcal{L}_{\mathbb{F}}^p([0,T]\times \Omega;E)$, and let $\psi^n$ be the martingale selection of $\Psi$ associated with $(x^n,g^n)$ for each $n$. By applying Doob's maximal inequality we have
\begin{eqnarray*}
\left\Vert \psi^m-\psi^n\right\Vert^p_{\mathscr{M}^{c,p}(E)}&=\quad & \mathbb{E}\left(\sup_{T\geq t\geq 0}\left\Vert x^m-x^n+\int_0^t\left(g^m_u-g^n_u\right) dW_u\right\Vert_E^p\right)\\
&\lesssim_{p\;\;} & \left\Vert x^m-x^n\right\Vert_E^p+\mathbb{E}\left(\left\Vert \int_0^T\left(g^m_u-g^n_u\right) dW_u\right\Vert_E^p\right)\\
&\eqsim_{p,E} & \left\Vert x^m-x^n\right\Vert_E^p+\left\Vert g^m-g^n\right\Vert^p_{L^p\left(\Omega;\gamma(0,T;E)\right)},
\end{eqnarray*}
for any natural numbers $n$ and $m$. Thus, $(\psi^n)_{n\geq 1}$ is a Cauchy sequence in $\mathscr{M}^{c,p}(E)$ that converges to $\psi$, and we have $\psi\in\mathscr{M}^p_{sel}(\Psi)$ by Lemma \ref{lem401}.\\Using the Burkholder–Davis–Gundy inequalities, we have, for every natural number $n$,
\begin{align*}
 \mathbb{E}\left(\sup_{T\geq t\geq 0}\left\Vert x+\int_0^tg_udW_u-\psi_t\right\Vert^p\right)
&\leq  2^{p-1}\left(\left\Vert x+\int_0^. g_udW_u-\psi^n\right\Vert_{\mathscr{M}^{c,p}(E)}^p+\left\Vert \psi^n-\psi\right\Vert_{\mathscr{M}^{c,p}(E)}^p \right)\\
& \lesssim_{p} \left\Vert x-x^n\right\Vert_E^p+\left\Vert\int_0^. g_udW_u-\int_0^. g^n_udW_u\right\Vert_{\mathscr{M}^{c,p}(E)}^p+\left\Vert \psi^n-\psi\right\Vert_{\mathscr{M}^{c,p}(E)}^p\\
&\eqsim_{p,E} \left\Vert x-x^n\right\Vert_E^p+\left\Vert g^n-g\right\Vert^p_{L^p\left(\Omega;\gamma(0,T;E)\right)}+\left\Vert \psi^n-\psi\right\Vert_{\mathscr{M}^{c,p}(E)}^p.
\end{align*}
Therefore, $(x,g)$ will be the element of $\mathfrak{S}_{\Psi}$ associated with $\psi$.\bigskip\\
Let $(x^\psi, f^\psi)$ be in $\mathfrak{S}_{\Psi}$. Since $\psi$ is assumed to have continuous trajectories, the Burkholder-Davis-Gundy inequalities once again lead us to obtain
\begin{eqnarray*}
\left\Vert f^{\psi}\right\Vert^p_{L^p(\Omega;\gamma(0,T;E))}&\;\,\eqsim_{p,E}& \left\Vert \psi-x^{\psi}\right\Vert_{\mathscr{M}^{c,p}(E)}^p\\
&\lesssim_{p}&  \left\Vert x^{\psi}\right\Vert_E^p+\left\Vert \psi\right\Vert_{\mathscr{M}^{c,p}(E)}^p\\
&\lesssim_{p}&	\left\Vert x^{\psi}\right\Vert_E^p+\left\Vert \upsilon\right\Vert_{L^{p}(\Omega)}^p.
\end{eqnarray*}
Besides this, the contraction property asserts that
$$\left\Vert x^{\psi}\right\Vert_E^p\leq \left\Vert \upsilon\right\Vert^p_{L^p(\Omega)}.$$
Thus, the boundedness of $\mathfrak{S}_{\Psi}$ has been established. It remains to note that the nature of the images taken by $\Psi$ ensures the convexity of the set $\mathscr{M}_{sel}^p(\Psi)$, which implies that $\mathfrak{S}_{\Psi}$ is convex.
\end{proof}
\begin{thm}\label{thm401}
If $\Psi$ is a set-valued martingale, then
\begin{align*}\label{ineq11}
\exists \mathfrak{S}\in \mathscr{K}_{cwcmpt}(E\times\mathcal{L}^p_{\mathbb{F}}\left([0,T]\times \Omega;E\right)),\; \forall t\in[0,T],\;  \Psi_t=\int_{[0,t]}^{\mathscr{R}}\mathfrak{S}\,  dW_u\quad a.s.\;\omega\in\Omega.
\end{align*}
\end{thm}
\begin{proof}
By virtue of Lemma \ref{lem403} we derive that the subset $\mathfrak{S}_{\Psi}$ defined by (\ref{ineq10}) belongs to the set $\mathscr{K}_{cwcmpt}(E\times\mathcal{L}^p_{\mathbb{F}}\left([0,T]\times \Omega;E\right))$. It follows according to Lemma \ref{lem402} that for every $t\in [0,T]$ we have
$$\overline{dec}_{\mathscr{F}_t}\left(\mathbb{I}_t(\mathfrak{S}_{\Psi})\right)=\overline{dec}_{\mathscr{F}_t}\left(\pi_t\left(\mathscr{M}^p_{sel}\left(\Psi\right)\right)\right)=S^p_{\mathscr{F}_t}\left(\Psi_t\right).$$
Therefore, the aim is obtained according to Corollary $1.2$ in \cite{Hiai}.
\end{proof}
\begin{rk}
By virtue of the tower property of conditional expectation, we observe that the set-valued stochastic process $\left(\int_{[0,t]}^{\mathscr{R}}\tilde{\mathfrak{S}},dW_u\right)_{T \geq t \geq 0}$ is a submartingale, in the sense of Definition 3.1 in \cite{jinping0},  for every nonempty bounded convex subset $\tilde{\mathfrak{S}}$ of the space $E \times \mathcal{L}^p_{\mathbb{F}}\left([0,T] \times \Omega; E\right)$.
\end{rk}
\noindent This prompts us to define the following set:
\begin{eqnarray}\label{assertion9}
	\mathbb{M}_{E\times\mathcal{L}^p_{\mathbb{F}}}:=\left\lbrace \tilde{\mathfrak{S}}\in \mathscr{K}_{cwcmpt}(E\times\mathcal{L}^p_{\mathbb{F}}\left([0,T]\times \Omega;E\right)):\, \left(\int_{[0,t]}^{\mathscr{R}}\tilde{\mathfrak{S}}\,dW_u\right)_{T\geq t\geq 0}\;\textrm{ is a martingale}\right\rbrace.	\end{eqnarray}
 
\subsection{Casting Representation}
\begin{cor}
If $\Psi$ is a set-valued martingale and $\left(\Omega, \mathscr{F}, \mathbb{P}\right)$ is separable, then there exists a sequence $\left((x^n, f^n)\right)_{n \geq 1}$ in $E \times \mathcal{L}^p_{\mathbb{F}}\left([0,T] \times \Omega; E\right)$ such that for every $t \in [0,T]$:
\begin{eqnarray}
\Psi_t &=& cl_E\left\lbrace x^n+\int_0^tf_u^ndW_u:n\geq 1\right\rbrace\quad a.s.\;\omega\in\Omega;\label{ineq12}\\
S^p_{\mathscr{F}_t}(\Psi_t)&=&\overline{dec}_{\mathscr{F}_t}\left\lbrace x^n+\int_0^tf^n_udW_u:n\geq 1\right\rbrace\label{ineq13}.
\end{eqnarray}
\end{cor}
\begin{proof}
By considering the set $\mathfrak{S}_{\Psi}$ defined in (\ref{ineq10}), we deduce, according to Theorem \ref{thm306}, that there exists a sequence $\left((x^n,f^n)\right)_{n\geq 1}$ contained in $\mathfrak{S}_{\Psi}$ such that for every $t\in [0,T]$:
$$\left(\int_{[0,t]}^{\mathscr{R}}\mathfrak{S}_{\Psi}\,  dW_u\right)(\omega)=cl_E\left\lbrace \left(x^n+\int_0^tf^n_udW_u\right)(\omega):n\geq 1\right\rbrace\quad a.s.\;\omega\in\Omega.$$
Thus, (\ref{ineq12}) follows by virtue of  Theorem \ref{thm401}.\\
To obtain (\ref{ineq13}), it suffices to proceed similarly to the final part in the proof of Theorem \ref{thm306}.
\end{proof}
\subsection{Hukuhara Difference of Set-valued Martingales}
We will explicitly use the set-valued conditional expectation. In such cases, the expression for its subtrajectory integrals must be specified due to the nature of the range of the set-valued functions involved.\\
Let $\Phi:\Omega\longrightarrow \mathscr{K}_{cwcmpt}(E)$ be an $\mathscr{F}_T$-measurable $p$-integrably bounded set-valued function.
\\Let $t\in [0,T]$. Let $(h^n)_{n\geq 1}$ be a sequence contained in $S^p_{\mathscr{F}_T}\left(\Phi\right)$ such that $\left(\mathbb{E}^{\mathscr{F}_t}(h^n)\right)_{n\geq 1}$ converging in $L^p(\Omega,\mathscr{F}_t;E)$. Since $(h^n)_{n\geq 1}$ is bounded in the reflexive space $L^p(\Omega,\mathscr{F}_T;E)$, hence it possesses a weak convergent subsequence $(h^{\varphi(n)})_{n\geq 1}$ to $h$ in $L^p(\Omega,\mathscr{F}_T;E)$. It follows by contractiveness of the conditional expectation that $\left(\mathbb{E}^{\mathscr{F}_t}(h^{\varphi(n)})\right)_{n\geq 1}$ is weakly convergent to $\mathbb{E}^{\mathscr{F}_t}(h)$ in $L^p(\Omega,\mathscr{F}_t;E)$.
Thus,
\begin{eqnarray}\label{assertion2}
\forall t\in [0,T],\quad S^p_{\mathscr{F}_t}\left(\mathbb{E}^{\mathscr{F}_t}\left(\Phi\right)\right)=\left\lbrace \mathbb{E}^{\mathscr{F}_t}\left(h\right):h\in S^p_{\mathscr{F}_T}\left(\Phi\right)\right\rbrace.
\end{eqnarray}
The following lemma leads us to justify the measurability of the Hukuhara difference of two set-valued functions and derive additional results. Prior to this, we recommend stating that the Hukuhara difference of two set-valued functions into $\mathscr{K}_{cwcmpt}(E)$ exists if, for every element in their common domain, the Hukuhara difference of their images exists.
\begin{lem}\label{lem405}
Let $\Phi^{(1)},\Phi^{(2)}:\left(\Omega,\mathscr{F}_{T}\right)\longrightarrow \mathscr{K}_{cwcmpt}(E)$ be two measurable set-valued functions.
\\If\, $\Phi^{(1)}\circleddash\Phi^{(2)}$ exists, then it is $\mathscr{F}_T$-measurable. Furthermore, if in addition $\Phi^{(1)}$ and $\Phi^{(2)}$ are $p$-integrably bounded, $\Phi^{(1)}\circleddash\Phi^{(2)}$ is $p$-integrably bounded, and we have $$\mathbb{E}^{\mathscr{F}_t}\left(\Phi^{(1)}\circleddash\Phi^{(2)}\right)=\mathbb{E}^{\mathscr{F}_t}\left(\Phi^{(1)}\right)\circleddash \mathbb{E}^{\mathscr{F}_t}\left(\Phi^{(2)}\right)\;a.s.\;\omega\in\Omega,\;t\in [0,T].$$
\end{lem}
\begin{proof}
Let $(e^*_n)_{n\geq 1}$ be a sequence from $E^*\setminus\left\lbrace 0_{E^*}\right\rbrace$ dense in $E^*$. Let $(\phi_m^{(1)})_{m\geq 1}$ (resp. $(\phi_m^{(2)})_{m\geq 1}$) be a sequence of measurable functions representing $\Phi^{(1)}$ (resp. $\Phi^{(2)}$) and let's put for every $n$  and every $\omega\in \Omega$:
\begin{eqnarray*}
h_n(\omega)&:=& \sup_{m\geq 1}\langle\phi_m^{(1)}(\omega),e^*_n\rangle-\sup_{m\geq 1}\langle\phi_m^{(2)}(\omega),e^*_n\rangle,\\
\Lambda_n(\omega)&:=& (-\infty,h_n(\omega)],\\
\Gamma_n(\omega)&:=&\left\lbrace e\in E:\langle e,e^*_n\rangle\,\leq h_n(\omega)\right\rbrace.
\end{eqnarray*}
Since for every real number $r$ and every $\omega\in \Omega$ we have: $d(r,\Lambda_n(\omega))=\mathds{1}_{\left\lbrace h_n(\omega)<r\right\rbrace}(r-h_n(\omega))$, thus $\Lambda_n$ is an $\mathscr{F}_T$-measurable set-valued function according to Theorem $3.5$ in \cite{Him}.\\
Let $\mathscr{O}$ be a nonempty open subset of $E$. Since $E$ is a separable metric space, then $\mathscr{O}$ is a countable union of balls $\left(\mathcal{B}(x_j,r_j)\right)_{j\geq 1}$. Subsequently, $\mathscr{O}$ is also union of the closed balls $$\overline{\mathcal{B}}_{j,m}:=\overline{\mathcal{B}}(x_j,(1-\frac{1}{m})r_j);\quad j\geq 1,\;m\geq 1.$$
For a fixed pair $(j,m)$ of natural numbers, let $(z_k)_{k\geq 1}$ be a converging sequence from $e^*_n(\overline{\mathcal{B}}_{j,m})$ and let $(y_k)_{k\geq 1}$ be from $\overline{\mathcal{B}}_{j,m}$ such that $e^*_n(y_k)=z_k$ for every $k$. According to the reflexivity of the space $E$ and the weak-closedness of $\overline{\mathcal{B}}_{j,m}$ we infer the existence of a subsequence $(y_{\varphi(k)})_{k\geq 1}$ weakly converging to an element of $ \overline{\mathcal{B}}_{j,m}$. Thus, $e^*_n(\overline{\mathcal{B}}_{j,m})$ is a closed set.\\Since $\Omega$ has a complete finite measure, the third property in [\cite{Him}, Theorem $3.5$ ] yields that:
\begin{eqnarray}\label{ineq21}
\forall (j,m),\quad\left\lbrace \omega\in \Omega:\Lambda_n(\omega)\cap e^*_n(\overline{\mathcal{B}}_{j,m})\neq \emptyset\right\rbrace\in \mathscr{F}_T.
\end{eqnarray}
On the other hand we have:
\begin{eqnarray}\label{ineq22}
\left\lbrace \omega\in \Omega:\Gamma_n(\omega)\cap \mathscr{O}\neq \emptyset\right\rbrace =\bigcup_{j,m\geq 1}\left\lbrace \omega\in \Omega:\Lambda_n(\omega)\cap e^*_n(\overline{\mathcal{B}}_{j,m})\neq \emptyset\right\rbrace.
\end{eqnarray}
Then, by combining (\ref{ineq21}) and (\ref{ineq22}) we derive that $\Gamma_n$ is an $\mathscr{F}_T$-measurable set-valued function.
\\Proposition \ref{prop203} asserts that
\begin{eqnarray*}
\Phi^{(1)}(\omega)\circleddash \Phi^{(2)}(\omega)= \bigcap_{n\geq 1}\left\lbrace e\in E:\langle e,e^*_n\rangle\,\leq \sup_{a\in \Phi^{(1)}(\omega)}\langle a,e^*_n\rangle-\sup_{b\in \Phi^{(2)}(\omega)}\langle b,e^*_n\rangle\right\rbrace =\bigcap_{n\geq 1}\Gamma_n(\omega),
\end{eqnarray*}
for every $\omega\in \Omega$. On the other hand, we have
\begin{align*}
	\Big\{ (\omega, x) \in \Omega \times E : x \in \bigcap_{n \geq 1} \Gamma_n(\omega) \Big\} 
	= \bigcap_{n \geq 1} \Big\{ (\omega, x) \in \Omega \times E : x \in \Gamma_n(\omega) \Big\}
\end{align*}
Thus, by applying  [\cite{Him}, Theorem $3.3$, Theorem $3.5$-(iii)] we obtain that the countable intersection $\Phi^{(1)}\circleddash \Phi^{(2)}$ is $\mathscr{F}_T$-measurable. Moreover, $\Phi^{(1)}\circleddash \Phi^{(2)}$ is $p$-integrably bounded by virtue of (\ref{formula1}).\\
Let $t$ be in $[0,T]$. Using (\ref{assertion2}) we conclude that
\begin{align*}
& S^p_{\mathscr{F}_t}\left(\mathbb{E}^{\mathscr{F}_t}(\Phi^{(2)})+\mathbb{E}^{\mathscr{F}_t}(\Phi^{(1)}\circleddash\Phi^{(2)})\right)\\& = cl_{L^p(\Omega,\mathscr{F}_t;E)}\left(S^p_{\mathscr{F}_t}\left(\mathbb{E}^{\mathscr{F}_t}(\Phi^{(2)})\right)+S^p_{\mathscr{F}_t}\left(\mathbb{E}^{\mathscr{F}_t}(\Phi^{(1)}\circleddash\Phi^{(2)})\right)\right)\\
&= S^p_{\mathscr{F}_t}\left(\mathbb{E}^{\mathscr{F}_t}(\Phi^{(2)})\right)+S^p_{\mathscr{F}_t}\left(\mathbb{E}^{\mathscr{F}_t}(\Phi^{(1)}\circleddash\Phi^{(2)})\right)\\
&=\left\lbrace \mathbb{E}^{\mathscr{F}_t}\left(h\right):h\in \left(S^p_{\mathscr{F}_T}\left(\Phi_2\right)+S^p_{\mathscr{F}_T}\left(\Phi^{(1)}\circleddash\Phi^{(2)}\right)\right)\right\rbrace\\
&=\left\lbrace \mathbb{E}^{\mathscr{F}_t}\left(h\right):h\in S^p_{\mathscr{F}_T}(\Phi^{(2)}+\Phi^{(1)}\circleddash\Phi^{(2)})\right\rbrace.
\end{align*}
Therefore, this is sufficient to conclude the aim.
\end{proof}
Next, we state the following representation theorem.
\begin{thm}
Let $\Psi^1,\Psi^2:[0,T]\times \Omega\longrightarrow\mathscr{K}_c(E)$ be two $\mathbb{F}$-adapted set-valued process such that $\Psi_T^1,\Psi_T^2$ are $p$-integrably bounded. If $\Psi_T^1 \circleddash\Psi_T^2$ exists and $\Psi^1,\Psi^2$ are set-valued martingales, then there exists $\mathfrak{S}$ in $\mathbb{M}_{E\times\mathcal{L}^p_{\mathbb{F}}}$ such that for every $t\in [0,T]$,
$$\Psi_t^1=\Psi_t^2+\int_{[0,t]}^{\mathscr{R}}\mathfrak{S}\, dW_s\quad a.s.\;\omega\in\Omega,$$
where $\mathbb{M}_{E\times\mathcal{L}^p_{\mathbb{F}}}$ is defined by (\ref{assertion9}).
\end{thm}
\begin{proof}
Taking into account Lemma \ref{lem405}, we deduce that  $\Psi_T^1 \circleddash\Psi_T^2$ is $\mathscr{F}_T$-measurable and $p$-integrably bounded. Subsequently, Theorem \ref{thm401} implies the existence of $\mathfrak{S}\in\mathbb{M}_{E\times\mathcal{L}^p_{\mathbb{F}}}$ such that, for every $t\in [0,T]$, we have
$$\mathbb{E}^{\mathscr{F}_t}\left(\Psi_T^1 \circleddash\Psi_T^2\right)=\int_{[0,t]}^{\mathscr{R}}\mathfrak{S}\, dW_u,\;a.s.\;\omega\in\Omega.$$
Finally, Lemma \ref{lem405} can be applied once more to complete the proof.
\end{proof}
%\begin{rk}
%Assume that ...As long as $\Psi_T^{(j)}$ is $p$-integrably bounded, let's choose $\upsilon^{(j)}$ from $L^p\left(\Omega,\mathscr{F}_T;E\right)$ such that $\left\vert \Psi_T^{(j)}\right\vert_{\varrho}\leq \upsilon^{(j)}\;a.s.\;\omega\in\Omega$.
%Let $t\in [0,T]$ and let $C\in \mathscr{F}_t$. By virtue of Theorem $5.2$ in \cite{Hiai} we conclude that
%\begin{eqnarray*}
%	\mathbb{E}\left(\mathds{1}_C\left\vert \Psi_t^{(j)}\right\vert_{\varrho}\right)=\mathbb{E}\left(\left\vert \mathbb{E}^{\mathscr{F}_t}\left(\mathds{1}_C\Psi_T^{(j)}\right)\right\vert_{\varrho}\right)&\leq & \mathbb{E}\left(\left\vert \mathds{1}_C\Psi_T^{(j)}\right\vert_{\varrho}\right)\\
%	&=& \mathbb{E}\left(\mathbb{E}^{\mathscr{F}_t}\left(\left\vert \mathds{1}_C\Psi_T^{(j)}\right\vert_{\varrho}\right)\right)\\
%	&\leq & \mathbb{E}\left(\mathds{1}_C\mathbb{E}_{con}^{\mathscr{F}_t}\left(\upsilon^{(j)}\right)\right).
%\end{eqnarray*}
%Consequently, for every $t\in [0,T]$ we have $\left\vert \Psi_t^{(j)}\right\vert_{\varrho}\leq \sup_{u\in [0,T]}\mathbb{E}_{con}^{\mathscr{F}_u}\left(\upsilon^{(j)}\right)\;a.s.\;\omega\in\Omega.$\\
%Thus, the set-valued martingale closed by $\Psi_T^1 \circleddash\Psi_T^2$ satisfies (\ref{assertion}). 	
%\end{rk}
\section{Existence of Solution for some Set-Valued Backward Differential Equations in UMD Spaces}
Let $E$ be a separable UMD space and assume that $\mathbb{F}$ is the natural filtration generated by a standard Brownian motion $W=\left(W_t\right)_{t\in [0,T]}$ augmented by all $\mathbb{P}$-null sets of $\mathscr{F}$.\\
We shall focus on the following form of Set-Valued Backward Stochastic Differential Equations (SVBSDE for short):
\begin{eqnarray}\label{eq501}
Y_t=\left(\xi+\int_t^TH_u du+\int_{[0,t]}^{\mathscr{R}}{Z}\,  dW_u \right)\circleddash \int_{[0,T]}^{\mathscr{R}}{Z}\,  dW_u\quad a.s.\;\omega\in\Omega,\quad t\in [0,T],
\end{eqnarray}
where $H:[0,T]\times \Omega\longrightarrow \mathscr{K}_c(E)$ be $\Sigma_{\mathbb{F}}$-measurable and $\xi:\Omega\longrightarrow \mathscr{K}_{cwcmpt}(E)$ be $\mathscr{F}_T$-measurable.

Let us now introduce the definition of our SVBSDE.
\begin{defi}
A pair $(Y,{Z})$ is called solution for the SVBSDE (\ref{eq501}) if it satisfies the following properties:
\begin{itemize}
\item[(a)]$Y$ is an adapted set-valued function with convex weakly compact valued;
\item[(b)]${Z}$ belongs to $\mathbb{M}_{E\times\mathcal{L}^p_{\mathbb{F}}}$.
\end{itemize}
\end{defi}

As a natural cognitive interaction with the formulation of this equation, the relationship between the case of set-valued functions and the simple case of functions must be clarified. For this purpose, let $h$ be a measurable and adapted mapping defining an element of $L^p_{\mathbb{F}}\left(\Omega;\gamma(0,T;E)\right)$, and let $\mathfrak{T}$ be an element of $L^p(\Omega,\mathscr{F}_T;E)$. For $H=\left\lbrace f\right\rbrace$ and $\xi=\left\lbrace \mathfrak{T}\right\rbrace$, the study made in the paper \cite{BSEE.B} assures us that there is a pair $(y,z)$ of adapted processes defining elements in $L^p_{\mathbb{F}}\left(\Omega;\gamma(0,T;E)\right)$ such that for every $t\in [0,T]$ we have
$$y_t=\mathfrak{T}+\int_t^Tf_u du-\int_t^Tz_u dW_u\quad\textrm{in}\;L^p(\Omega;E).$$
Thus, for every $t\in [0,T]$
$$\left\lbrace y_t\right\rbrace=\left(\xi+\int_t^T\left\lbrace f_u\right\rbrace du+\int_{[0,t]}^{\mathscr{R}}\left\lbrace (0_E, z)\right\rbrace dW_u\right) \circleddash \int_{[0,T]}^{\mathscr{R}}\left\lbrace (0_E, z)\right\rbrace dW_u\quad a.s.\;\omega\in\Omega.$$
Equivalently, the pair $(Y,Z)=(\left\lbrace y\right\rbrace,\left\lbrace (0_E, z)\right\rbrace)$ solves the problem (\ref{eq501}) in this special case.
\\Here, it should be noted that the pair $(Y,\tilde{Z})=(\left\lbrace y\right\rbrace,\left\lbrace (x, z)\right\rbrace)$ solves also (\ref{eq501}) for every $x\in E$.
\bigskip

Next, let us consider the following assumptions:
\begin{itemize}
\item[$(\mathscr{H}_1)$] The set $\Lambda_{p,H}^{0,T}$ defined by (\ref{formula2}) is nonempty and satisfies assumption (\ref{assertion6}).
\item[$(\mathscr{H}_2)$]$\xi$ is $\mathscr{F}_T$-measurable and $p$-integrably bounded.
\end{itemize}
\begin{thm}\label{thm501}
If the hypotheses $(\mathscr{H}_1)-(\mathscr{H}_2)$ are satisfied, then there exists a pair $(Y,Z)$ solving the SVBSDE (\ref{eq501}), where $Y$ is unique up to modifications.
\end{thm}
\begin{proof}
Clearly, the process defined as follows
\begin{eqnarray}\label{ineq17}
\Psi_t:=\mathbb{E}^{\mathscr{F}_t}\left(\xi+\int_0^TH_u du\right),
\end{eqnarray}
is a set-valued martingale satisfying for every $t\in [0,T]$
$$\left\vert \Psi_t\right\vert_{\varrho}\leq \sup_{u\in [0,T]}\mathbb{E}_{con}^{\mathscr{F}_u}\left(\left\vert\xi\right\vert_{\varrho}\right)+\sup_{s\in [0,T]}\mathbb{E}_{con}^{\mathscr{F}_s}\left(\left\vert \int_0^TH_u du\right\vert_{\varrho}\right)\;a.s.\;\omega\in\Omega.$$
Thus, by virtue of Theorem \ref{thm401} we derive the existence of a nonempty convex closed bounded subset $Z$ contained in $E\times\mathcal{L}^p_{\mathbb{F}}([0,T]\times \Omega;E)$ such that for every $t\in [0,T]$
\begin{eqnarray}\label{ineq18}
\Psi_t=\int_{[0,t]}^{\mathscr{R}} Z\,  dW_u\quad a.s.\;\omega\in\Omega.
\end{eqnarray}
Let us put for every $t\in [0,T]$
$$Y_t:=\mathbb{E}^{\mathscr{F}_t}\left(\xi+\int_t^TH_u du\right).$$
Taking into account Theorem \ref{thm311}, we obtain
\begin{eqnarray}\label{ineq19}
\Psi_t=\mathbb{E}^{\mathscr{F}_t}\left(\xi+\int_0^t H_u du+\int_t^TH_u du\right)= Y_t+\int_0^t H_u du\quad a.s.\;\omega\in\Omega.
\end{eqnarray}
It follows by considering inequalities (\ref{ineq18})-(\ref{ineq19}), that
\begin{eqnarray}\label{ineq20}
\Psi_T+\int_{[0,t]}^{\mathscr{R}} Z\,  dW_u &=& \Psi_t+\int_{[0,T]}^{\mathscr{R}} Z\,  dW_u \nonumber\\
&=& Y_t+\int_0^t H_u du+\int_{[0,T]}^{\mathscr{R}} Z\,  dW_u\quad a.s.\;\omega\in\Omega.
\end{eqnarray}
So, by using (\ref{ineq17}) and substituting it into (\ref{ineq20}), we obtain
\begin{eqnarray*}
\xi+\int_0^T H_u du +\int_{[0,t]}^{\mathscr{R}} Z\,  dW_u=Y_t+\int_0^tH_u du+\int_{[0,T]}^{\mathscr{R}} Z\,  dW_u\quad a.s.\;\omega\in\Omega.
\end{eqnarray*}
Next, by virtue of the cancellation law lemma we conclude that
\begin{eqnarray}\label{eq502}
\xi+\int_t^TH_u du +\int_{[0,t]}^{\mathscr{R}} Z\,  dW_u=Y_t+\int_{[0,T]}^{\mathscr{R}} Z\, dW_u\quad a.s.\;\omega\in\Omega.
\end{eqnarray}
Therefore, by Definition \ref{defi201} we get the existence of a solution.\\
Let $(\tilde{Y},\tilde{Z})$ be a pair of another solution to (\ref{eq501}) and let $t\in [0,T]$. From (\ref{eq502}) we conclude that
\begin{eqnarray}\label{eq503}
\left(Y_t+\int_{[0,T]}^{\mathscr{R}} Z\, dW_u\right) \circleddash \int_{[0,t]}^{\mathscr{R}} Z\, dW_u=\left(\tilde{Y}_t+\int_{[0,T]}^{\mathscr{R}} \tilde{Z}\, dW_u\right) \circleddash \int_{[0,t]}^{\mathscr{R}} \tilde{Z}\, dW_u\quad a.s.\;\omega\in\Omega.
\end{eqnarray}
Thus, by taking into account Lemma \ref{lem405} and the adaptness of $Y^{(j)},j=1,2$, we derive after taken conditional expectation with respect to $\mathscr{F}_t$ that $Y_t^{(1)}=Y_t^{(2)}\;a.s.\;\omega\in\Omega.$
\end{proof}
\begin{rk}
In addition to what has been achieved in the point-valued case, one can naturally inquire about the second component in a multivalued case, specifically concerning two pairs of solutions $(Y,Z^{(1)})$ and $(Y,Z^{(2)})$ to (\ref{eq501}). Assume that
\begin{eqnarray}\label{eq504}
\left\{
	\begin{array}{lll}
		 (i) &
\displaystyle\int_{\left\lbrace 0\right\rbrace}^{\mathscr{R}} Z^{(1)}\, dW_u	=\int_{\left\lbrace 0\right\rbrace}^{\mathscr{R}} Z^{(2)}\, dW_u\quad a.s.\;\omega\in\Omega
	 ;\bigskip
	
\\ (ii) & 	\text{the Hukuhara difference } Z^{(1)}\circleddash Z^{(2)} \text{ exists}.
		
	\end{array}
	\right.
\end{eqnarray}
By combining (\ref{eq503}), (\ref{eq504})-(i) and using the cancellation law lemma we derive that
$$\int_{[0,T]}^{\mathscr{R}} Z^{(1)}\, dW_u=\int_{[0,T]}^{\mathscr{R}} Z^{(2)}\, dW_u\quad a.s.\;\omega\in\Omega.$$
Let $(x,f)$ be from $E\times\mathcal{L}_{\mathbb{F}}^p([0,T]\times \Omega;E)$ such that $(x,f)+Z^{(2)}\subset Z^{(1)}$. One can see obviously that
$$\mathbb{I}_T\left((x,f)\right)+S^p_{\mathscr{F}_T}\left(\int_{[0,T]}^{\mathscr{R}}Z^{(2)}\,dW_u\right)\subset S^p_{\mathscr{F}_T}\left(\int_{[0,T]}^{\mathscr{R}}Z^{(1)}\,dW_u\right).$$
Thus, $x+\int_0^Tf_udW_u=0_{L^p(\Omega,\mathscr{F}_T;E)}$ and $(x,f)=(0_E,0_{\mathcal{L}_{\mathbb{F}}^p([0,T]\times \Omega;E)})$. Therefore, $Z^{(1)}$ and $Z^{(2)}$ are equal as a subsets of $E\times\mathcal{L}_{\mathbb{F}}^p([0,T]\times \Omega;E)$.
\end{rk}
\par There was considerable interest concerning stochastic integration of martingale type 2 spaces, which influenced some authors to invest in multivalued analysis (see \cite{jinping}, for example).
\\Here, we recall that a UMD space
has martingale type $2$ if and only if it has type $2$ (cf. \cite{neerven3}).
\bigskip\\Now, to derive the next result, we consider the following hypothesis.
\begin{itemize}
\item[$(\mathscr{H}_3)$]The set-valued function $H$ has the following property:
$$\displaystyle\int_{\Omega}\left\Vert \vert H(.,\omega)\vert_{\varrho}\right\Vert_{L^2(0,T)}^2\mathbb{P}(d\omega)<\infty.$$
\end{itemize}
\begin{cor}
Let $p=2$ and assume that the hypotheses $(\mathscr{H}_2)-(\mathscr{H}_3)$ are satisfied. If $E$ has martingale type $2$, then the set-valued backward stochastic differential equation (\ref{eq501}) admits a solution.
\end{cor}
\begin{proof}
Firstly, according to  Kuratowski and Ryll-Nardzewski selection theorem [\cite{Him}, Theorem $5.1$] we guarantee the existence of a measurable and adapted selection $h:[0,T]\times \Omega\longrightarrow E$ of the set-valued function $H$. It follows by virtue of $(\mathscr{H}_3)$ that $h\in L^2(\Omega;L^2(0,T;E))$. Thus, Theorem $9.2.10$ in \cite{neerven2} asserts that $h$ belongs to $L^2(\Omega;\gamma(0,T;E))$, and the set $\Lambda_{2,H}^{0,T}$ is nonempty.\\
Let $g\in \Lambda_{2,H}^{0,T}$ and let $A\in \mathscr{F}_T$. Once more, Theorem $9.2.10$ in \cite{neerven2} leads us to obtain 
$$(\varsigma_{2,E}^{\gamma})^2\int_A\left\Vert \vert H(.,\omega)\vert_{\varrho}\right\Vert_{L^2(0,T)}^2\mathbb{P}(d\omega)\geq (\varsigma_{2,E}^{\gamma})^2\int_A\left\Vert g(.,\omega)\right\Vert_{L^2(0,T;E)}^2\mathbb{P}(d\omega)\geq \int_A\left\Vert g(.,\omega)\right\Vert_{\gamma(0,T;E)}^2\mathbb{P}(d\omega),$$
where $\varsigma_{2,E}^{\gamma}$ is the Gaussian type $2$ constant (see \cite{neerven2}, Definition 7.1.17). Then,
$$\forall g\in \Lambda_{2,H}^{0,T},\;\left\Vert g(.,\omega)\right\Vert_{\gamma(0,T;E)}\leq \varsigma_{2,E}^{\gamma}\left\Vert \vert H(.,\omega)\vert_{\varrho}\right\Vert_{L^2(0,T)}\;a.s.\;\omega\in\Omega.$$
Finally, to accomplish the proof it suffices to apply Theorem \ref{thm501}.
\end{proof}
\section*{Conclusion}
In this paper, we establish a new tool that links multivalued analysis and its applications to stochastic integration in UMD spaces. We present several results and properties that, to the best of our knowledge, have not been previously explored in this context.
\\The revised set-valued stochastic integral enabled us to derive a martingale representation, which contributes to the solution of the SVBSDE (\ref{eq501}). However, we emphasize that solving this equation with $\tilde{H}(\cdot, Y, {Z})$  instead of $H$ presents significant challenges.\\
On the other hand, the subset $\Lambda_{p,H}^{0,T}$ has played a crucial role in defining the new concept "$\gamma$-set-valued stochastic integral of a set-valued function $H$". This subset can provide the framework for defining the set-valued stochastic integral of $H$ with respect to a standard Brownian motion, extending the study in [\cite{kis2012}, Section 3], and enabling the analysis of certain types of stochastic inclusions.

%\noindent On the other hand, although the revised set-valued stochastic integral is a component in the SVBSDE (\ref{eq501}) with a proven solution via the martingale representation, it doesn't seem easy to resolve the equation with $H(\cdot, Y, \mathcal{Z})$ instead of $H$.
% Additionally, the Lebesgue space $L^r(\mathcal{S}):=L^r(\mathcal{S}, \mathscr{G}, \nu)$ is a separable UMD space for every $r \in (1,\infty)$ and $\sigma$-finite measure space $(\mathcal{S}, \mathscr{G}, \nu)$ such that $\mathscr{G}$ is a $\nu$-countably generated $\sigma$-algebra. Hence, a special study can be carried out by virtue of the isomorphism from the $\gamma$-radonifying space $\gamma(0,T;L^r(\mathcal{S}))$ onto the Bochner space $L^r(\mathcal{S};L^2([0,T]))$ (see Proposition 9.3.2 in \cite{neerven2}).

\end{document}